\documentclass[a4paper]{amsart}

\usepackage[usenames,dvipsnames,table]{xcolor}
\usepackage{amsmath}
\usepackage{amssymb}
\usepackage{amsfonts}
\usepackage{esint}
\usepackage{subfig}
\usepackage{graphicx}
\usepackage{amscd}
\usepackage{hyperref}
\usepackage{comment}
\usepackage{ifthen}
\usepackage[numbers,sort&compress]{natbib}

\newtheorem{Defin}{Definition}
\newtheorem{Theo}{Theorem}
\newtheorem{Lemma}{Lemma} 
\newtheorem{Prop}{Proposition}
\newtheorem{Conjecture}{Conjecture}

\newtheorem{remark}{Remark}
\DeclareMathOperator{\Div}{\nabla \cdot}
\newfont{\NUMBERS}{msbm8 scaled\magstep1}
\newcommand{\REAL}{\mbox{\NUMBERS R}}
\newcommand{\Cac}{\mathcal{C}^{\delta}(\bar{\Omega})}
\newcommand{\Hac}{\mathcal{C}^{1,\delta}(\bar{\Omega})}
\newcommand{\cam}[2]{\mathcal{L}^{2,#1}(#2)}
\newcommand{\Supp}{\ensuremath{\operatorname{supp}}}
\newcommand{\mor}[2]{L^{2,#1}(#2)}
\newcommand{\Pot}{\mathcal{U}}
\newcommand{\Flux}{\mathcal{Q}}
\newcommand{\tdens}{\mu}
\newcommand{\Kdelta}{K_{\delta}}
\newcommand{\Kmu}{K_{\tdens}}
\newcommand{\dx}{\, dx}
\newcommand{\dy}{\, dy}
\newcommand{\ds}{\, ds}
\newcommand{\Lyap}{\mathcal{L}}
\newcommand{\Span}{\ensuremath{\operatorname{span}}}
\newcommand{\var}{\ensuremath{\operatorname{var}}}
\newcommand{\Tsymb}{{\mathcal T}}
\newcommand{\Triang}[1][]
{
  \ifthenelse{\equal{#1}{}}
  {\Tsymb}
  {\Tsymb_{#1}}
}
\newcommand{\PONE}{\ensuremath{\mathcal{P}_{1}}}
\newcommand{\PTWO}{\ensuremath{\mathcal{P}_{2}}}
\newcommand{\PZERO}{\ensuremath{\mathcal{P}_{0}}}
\newcommand{\Deltat}{\Delta t_{j}}
\newcommand{\BallSymb}{B}

\newcommand{\Ball}[2][]
{
  \ifthenelse{\equal{#2}{}}
  {\BallSymb_{#1}}
  {\BallSymb(#2,#1)}
}
\newcommand{\Avgint}[3][]
{
  \ifthenelse{\equal{#1}{}}
  {({#3})_{#2}}
  {(#3)_{#1,#2}}
}

\begin{document}

\markboth{E. Facca et al.}{
  Towards a stationary Monge-Kantorovich dynamics
}

\title[Towards a stationary Monge-Kantorovich dynamics]{
  Towards a stationary Monge-Kantorovich dynamics: the Physarum
  Polycephalum \\experience
}

\author{ENRICO FACCA, FRANCO CARDIN and MARIO PUTTI}
\address{Department of Mathematics, University of Padua, \\
  via Trieste 62, Padova, Italy\\
  \{facca,cardin,putti\}@math.unipd.it}

\maketitle

\begin{abstract}
In this work we study and expand a model describing the dynamics of
a unicellular slime mold, Physarum Polycephalum (PP), which was
proposed to simulate the ability of PP to find the shortest path
connecting two food sources in a maze. 
The original model describes the dynamics of the 
slime mold on a finite dimensional planar graph using a pipe-flow
analogy whereby mass transfer occurs because of pressure
differences with a conductivity coefficient that varies with the
flow intensity. 
We propose an extension of this model that abandons the graph
structure and moves to a continuous domain.
The new model, that couples an elliptic equation enforcing
PP density balance with an ODE governing the
dynamics of the flow of information along the PP body, is analyzed by
recasting it into an infinite-dimensional dynamical system. 
We are able to show well-posedness of the proposed model for
sufficiently small times under the hypothesis of H\"older continuous
diffusion coefficients and essentially bounded forcing functions,
which play the role of food sources. 
Numerical evidence, shows that the model is capable of describing the
slime mold dynamics also for large times, accurately reproducing the
PP behavior. 

A notable result related to the original model is that it is
equivalent to an optimal transportation problem over the graph as time
tends to infinity.  In our case, we can only conjecture that our
extension presents a time-asymptotic equilibrium.  This equilibrium
point is precisely the solution of the Monge-Kantorovich (MK)
equations at the basis of the PDE formulation of optimal
transportation problems.  Numerical results obtained with our
approach, which combines \PONE\ Finite Elements with forward Euler
time stepping, show that the approximate solution converges at large
times to an equilibrium configuration that well compares with the
numerical solution of the MK-equations.
\end{abstract}



\section{Introduction}

\begin{figure}
  \centerline{
    \includegraphics[width=0.5\columnwidth]{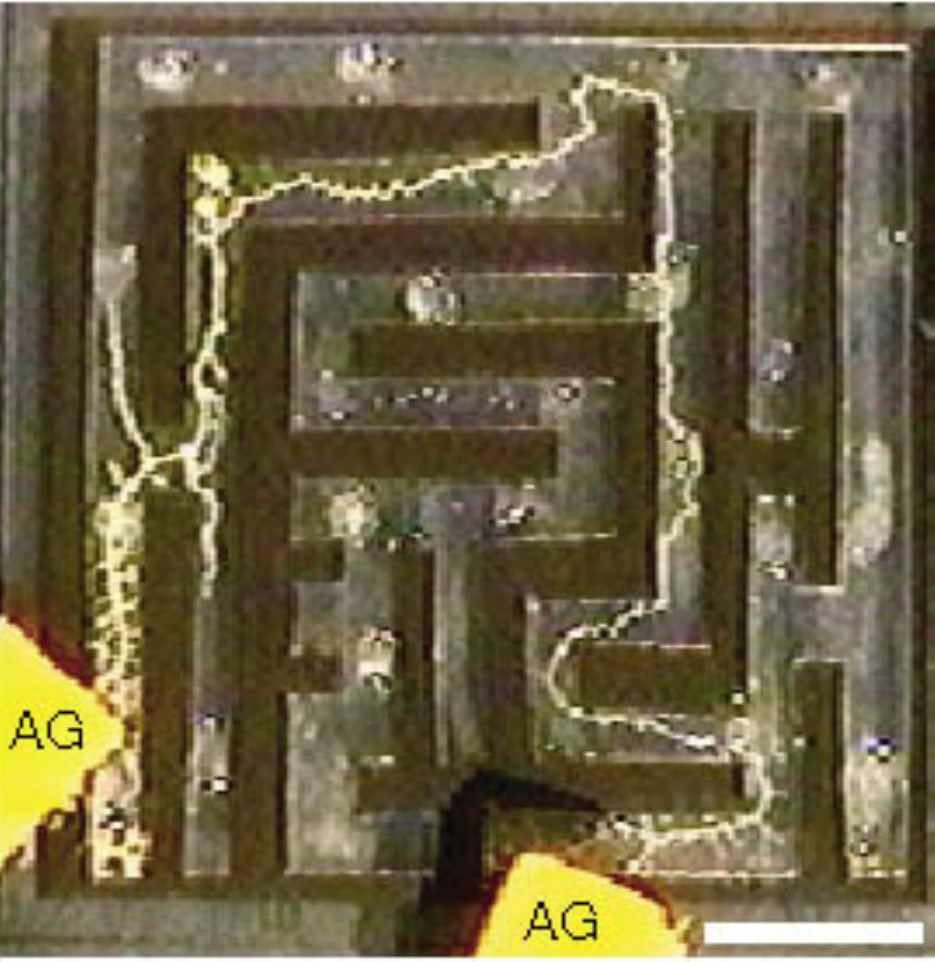}
  }
  \caption{\emph{Physarum Polycephalum} experimental maze solving
    (from~\citet{nakagaki:maze})} 
  \label{slime_maze}
\end{figure}

In a recent paper, \citet{tero:model} proposed a mathematical model
governing the dynamics of a unicellular slime mold named
\emph{Physarum Polycephalum} (PP) that, on the basis of experimental
evidence~\citep{nakagaki:maze}, grows following the most efficient
network path between food sources.  This evidence is exemplified in
the picture shown in Fig.~\ref{slime_maze} that portrays the
experimental setup developed by~\citet{nakagaki:maze} suggesting that
Physarum slime, after colonization of the entire maze paths, evolves
along the shortest path connecting the two food sources.
P. Polycephalum abilities have been used for the experimental analysis
of transportation networks, with many researchers suggesting that this
slime mold is capable of identifying the optimal many-site connecting
transportation network, with applications to such systems as the
railroads of Tokyo and Spain~\citep{tokyo,physarum:machine}.  Many
further surprising properties of PP have been
experimentally identified, but in this work we are interested in
studying and extending the mathematical model proposed
by~\citet{tero:model} that governs the slime mold dynamics.

The original model of~\citet{tero:model} reads as follows. Given a
undirected planar graph $G=(E,V)$, with positive edge length 
$\{L_e\}_{e\in E}$ find the edge function $D_e$ and the vertex function
$p_v$ that satisfy: 
\begin {subequations}
  \label{sys:finite}
  \begin{align}
    &\sum _{e \in \sigma (v)} Q_{e}(t)= f_{v} 
    \quad \text{(``Balance law-Kirchhoff'')}&&
        \forall 
        v\in V, \label{sys:finite:cons}\\
    &Q_{e}(t)=D_{e}(t)\frac{(p_{u}(t)-p_{v}(t))}{L_{e}} 
    \quad\text{(``Fick-Poiseuille'')}&&
        \forall 
        e=(u,v)\in E, \label{sys:finite:ohm}\\
    &D'_{e}(t)=g\left(\left|Q_{e}(t)\right|\right) -D_{e}(t)  
    \quad\text{($D_{e}$ dynamics)}&&
        \forall 
        e=(u,v)\in E, \label{sys:finite:dyn}\\
    &D_e(0)=\hat{D}_e(0)>0  
    \quad\text{(initial data)}&&
        \forall
        e=(u,v)\in E,
  \end{align}
\end{subequations}
where $e=(u,v)$ denotes the edge of $G$ connecting vertices $u$ and
$v$, the vertex source function $f_{v}$ satisfies the compatibility
condition of an isolated system
$\sum_{v\in V}f_v=0$~\citep{bonifaci:physarum}, where $\sigma(v)$ is
the ``star'' of $v$, i.e., the set of edges having vertex $v$ in
common and $g:\REAL^+\mapsto\REAL^+$ is a non-decreasing function with
$g(0)=0$.  This model can be explained heuristically using a classical
hydraulic analogy, eventually motivating the above introduced terms
``balance law-Kirchhoff'' and ``Fick-Poiseuille''. We think of the
graph $G$ as representing the set of pipes where the flow of a fluid
driven by the vertex source function $f_v$ occurs.  Then, the first
equation~\eqref{sys:finite:cons} can be identified as the enforcement
of the fluid mass balance, while equation~\eqref{sys:finite:ohm} is
the momentum balance stating that the flux in each graph edge is
proportional to the discrete gradient of the vertex potential function
$p_{v}$~\eqref{sys:finite:ohm} via a conductance coefficient $D_e$
(inverse of a resistance).  Hydraulic resistance to flow is known to
be proportional to the pipe perimeter, and hence to its
diameter. Thus, the evolutive equation~\eqref{sys:finite:dyn}, which
forms the core of the model, asserts the intuitive behavior that to
optimally (with minimal energy loss) accommodate larger fluxes, the
pipe diameter must increase, although it needs to remain bounded.
From this observation it can be concluded that the function $g(x)$
must be nondecreasing. Moreover, to avoid unboundedness, the growth of
the hydraulic conductivity needs to be compensated by introducing the
balancing decay term $-D_{e}(t)$.  In~\citet{tero:model} several
numerical results using this model were presented. The most relevant
to our study are those where the vertex source function was
concentrated in the first (1) and last ($n$) vertex of the graph (the
entrance and exit of the fluid in the hydraulic analogy), i.e.:
\begin{equation*}
  f_v=
  \begin{cases}
    +1 \qquad & v=1 \\
    -1        & v=n \\
    0         & \mbox{otherwise} 
 \end{cases}
\end{equation*}
Using this setting, the numerical evidence shows that when $g(x)=x$
the conductivity $D_e$ at large times tends to localize (have a
nonzero support) on the edges of the shortest path between the two
external sources. This has been confirmed more recently
by~\citet{bonifaci:physarum}, who show that, for $t\rightarrow\infty$,
indeed the distribution of $D_e$ converges to the shortest
path. Moreover, the same authors prove that the above model is
equivalent to an optimal transport problem on the graph $G$ and can be
recasted as the problem of finding $Q=\{Q_{e}\}_{e \in E}$ such that:
\begin{align*}
  & \min_{Q\in\{Q_{e}\}_{e\in E}} J(Q), \qquad J(Q) := \sum _{e\in E} Q_{e}L_{e} \\
  &\mbox{s.t.:} \\
  & \sum _{e \in \sigma (v)} Q_{e}= f_{v} \qquad  \mbox{for all } v\in V.
\end{align*}
In fact, under some general assumptions on the graph structure, 
the solution of system~\eqref{sys:finite} converges to
a stationary solution $\bar{Q}$ that is also solution of the above
optimal transport problem in $G$.

In this work we generalize the model given in~\eqref{sys:finite} by
removing the graph structure and defining 
the problem on an open bounded domain $\Omega\subset\REAL^n$. We
restrict this study to the case of $g(x)=x$. Then, given a source
function $f:\Omega\to\REAL$, a continuous analogue
of~\eqref{sys:finite} tries to find the pair of functions
$(\tdens,u):[0,+\infty[\times\Omega\mapsto\REAL^{+}\times\REAL^n$ that
satisfies:
\begin{subequations}
  \label{sys:intro}
  \begin{align}
    &  -\Div \Big( \tdens(t,x)\nabla  u(t,x)\Big) = f(x) 
      \qquad \left(\int_{\Omega} f\dx=0\right) 
    &\label{sys:intro:div}\\[-0.3em]
    & \tdens'(t,x) = \tdens(t,x)\Big(|\nabla u(t,x)|-1\Big)
    \label{sys:intro:dyn}\\
    & \tdens(0,x) =\tdens_0(x)> 0\label{bound:cond:d}
  \end{align}
\end{subequations}
complemented by zero Neumann boundary conditions.  Here, $\tdens'$
indicates partial differentiation with respect to time, and
$\nabla =\nabla_x$.  
This generalization is intuitively
justified by comparing the different components of
models~\eqref{sys:finite} and~\eqref{sys:intro}.  In
fact,~\eqref{sys:intro:div} states the spatial balance of a 
(continuum) Fick-Poiseuille flux $q=-\tdens\nabla u$ with potential
function $u$, while~\eqref{sys:intro:dyn} introduces
the dynamics introduced in the original discrete model.  To address
the dynamics of PP in the maze, we need to reconcile the model with the
fact that some portions of the domain (the maze barriers in this case)
may hinder through-flow. This can be obtained by imposing the gradient
to be large where the flux must be small, thus forcing the
conductivity $\tdens$ to become small.  Thus,
equation~\eqref{sys:intro:dyn} is replaced by:
\begin{equation}
  \tdens'(t,x)=\tdens(t,x)\Big(|\nabla u(t,x)| - k(x)\Big) 
  \label{het:intro:dyn} 
\end{equation}
where $k(x)$ is a positive function describing the spatial pattern of
the resistance to flow, whereby large values of $k$ imply large energy
losses and hence large gradients of the potential $u$. Note that,
intuitively, since the flux $q$ must remain bounded, larger gradients
induce smaller conductivities and hence smaller fluxes. 

In the present paper, we analyze from an analytical and numerical
point of view the continuous model of slime-mold dynamics described
above, and we conjecture that, like its discrete counterpart, its
solution tends to an equilibrium point as time goes to infinity and
that this equilibrium point is the solution of the Monge-Kantorovich
optimal mass transport problem problem~\citep{evans2}.   
We first study existence and uniqueness of the solution
of~\eqref{sys:intro} for the case $k\equiv 1$, as the more general
case with heterogeneous $k(x)$ is a straight forward adaptation, at
least for smooth enough $k$.  
For $k\equiv 1$, if an equilibrium
point exists for $t\rightarrow+\infty$, then $\tdens'\rightarrow
0$. Hence,~\eqref{sys:intro:dyn} becomes a
constraint imposing that for, $\tdens$ strictly greater than or equal
to zero, the norm of the gradient of $u$ must be unitary.  This
observation is crucial to the development of our conjecture, which
reads as: 
\begin{Conjecture}
  The solution of~\eqref{sys:intro} tends as $t\rightarrow\infty$ to the
  solution of the following problem: find
  $(\tdens^*,u^*)\in (L^{\infty}(\Omega), Lip_1(\Omega))$, with
  $Lip_1(\Omega)$ the space of Lipschitz continuous functions with
  unit constant, such that:
  \begin{align}
    &-\Div(\tdens^*(x) \nabla u^*(x))=f^+(x) -f^-(x)=f(x) 
    &&
       \mbox{ in } \Omega \nonumber \\
    &|\nabla u^*(x)| \le 1   
    &&\text{ in } \Omega \label{evans} \\
    &|\nabla u^*(x)| =1    
    && \text{ where } \tdens^*(x)>0 \nonumber 
  \end{align}
\end{Conjecture}
These equations constitute the PDE-based formulation of the classical
optimal transport problem and are named the \emph{Monge-Kantorovich}
(MK) equations~\citep{shape,evans2,Ambrosio:2004}.  In the past few
years they have been the subject of a number of studies
that have shed light into regularity and integrability properties of
the \emph{optimal transport density} $\tdens^*$ in relation to the
regularity and integrability of the forcing function
$f$~\citep{DePascale:2004}, and its uniqueness~\citep{Feldman:2002}.

In this work we show the applicability of the proposed model to the
simulation of the dynamics of Physarum Polycephalum, but, most
notably, we try to point out to some theoretical and numerical
evidence in support of veridicity of the above conjecture.  From a
theoretical point of view, we first prove the local-in-time existence
and uniqueness of the solution pair $(\tdens,u)$ in H\"older spaces.
To this aim, we recast the problem in operatorial form and look for
the functions $u(t)$ and $\tdens(t)$ such that:
\begin{subequations}
  \begin{align}
    & \int_\Omega \tdens(t)\nabla u(t)\nabla \varphi(x)\dx =
      \int_\Omega f\varphi(x)\dx \qquad 
      \forall \varphi\in H^1(\Omega) \label{pre:pde}\\
    & \tdens'(t)=\Flux\left(\tdens(t)\right)-\tdens(t); \qquad
      \tdens(0)=\tdens_0 \label{pre:ic} 
  \end{align}
\end{subequations}
where $\Flux(\tdens)$ is the operator, associated with the weak
form~\eqref{pre:pde} of the elliptic PDE~\eqref{sys:intro:div}, that
maps $\tdens$ into $\tdens|\nabla u|$.  To simplify notation, we
consider the dependence on $x\in\Omega$ implicit in all the relevant
functions.  For a given $\tdens>0$, we denote by $\Pot(\tdens)$ the
unique weak solution of the elliptic PDE associated with $\tdens$.
Clearly, the solution $\tdens(t)$ remains bounded as long as the
left-hand-side remains negative, i.e.,
$|\nabla \Pot(\tdens(t))|\le 1$.  Standard regularity theory of
elliptic PDEs ensures that, for any source function
$f\in L^\infty(\Omega)$ the operator $\Flux(\tdens)$ is well-defined
(meaning problem~\eqref{sys:intro:div} is well-posed, i.e., the
associated bilinear form is coercive and continuous) for $\tdens(t,x)$
H\"older continuous and strictly greater than zero.  We are then able
to show that the operator $\Flux(\tdens)$ is locally Lipschitz
continuous, and we can then invoke Banach-Caccioppoli fixed-point
theorems to show local existence and uniqueness of the solution pair
$(\tdens,u)$. However, the fact that Lipschitz continuity is only
local in $\tdens$, which is a consequence of the need to maintain
coercivity of the diffusion equation~\eqref{sys:intro:div}, prevents
the extension of this result to larger times.  Nevertheless we are
able to identify a Lyapunov-candidate function, i.e., a function whose
time derivative along the flow trajectories (its Lie derivative) is
strictly negative and is equal to zero only when the it attains its
minimum.

Next, on the basis of these findings, we experiment numerically the
large time behavior of the proposed system~\eqref{sys:intro} together
with the extended dynamics~\eqref{het:intro:dyn}.  Given a
triangulation of the domain $\Omega$, the discrete model uses a
$\PONE$ Galerkin approach for the discretization of the potential $u$
and a $\PZERO$ Galerkin scheme for the discretization of the diffusion
coefficient $\tdens$. To avoid oscillations on the numerical gradient,
we use different approximation spaces for $u$ and $\tdens$ by defining
the $\PONE$ potential on a uniformly refined mesh, an approach
reminiscent of the inf-sup stable $\PONE$-iso-$\PTWO$/$\PONE$ Stokes
FEM spaces. This methodology leads to a coupled differential algebraic
system of equations that is solved by a simple forward Euler method.
The resulting numerical algorithm is applied for the solution of the
dynamics of P. Polycephalum on the maze.

The conjecture that problems~\eqref{sys:intro} and~\eqref{evans} are
asymptotically equivalent finds application to the numerical solution
of the OT equations as, from preliminary experimental evidence,
convergence to an asymptotic state seems smooth and fast even using
the low order standard \PONE/\PZERO\ Galerkin approximations described
above and on relatively coarse meshes.  The conjecture is supported
experimentally by comparing our numerical solution against the
numerical results presented in~\citet{prigozhin}, where solution
to~\eqref{evans} was obtained by means of RT0 mixed finite elements
and automatic mesh refinement at on a larger triangulation and with
much higher computational costs.
Eventually, from a numerical point of view, the above asymptotic
behavior seems rather robust under iterated mesh refinements.

\section{Main results}

In this Section we lay down the technical results about our slime mold 
dynamics gained in a rather sharp and convenient functional
environment. The main idea consists on the synthesis of
the proposed dynamics towards an infinite dimensional abstract ODE
that passes through a standard elliptic setting. We obtain a theorem
about local existence and uniqueness, but at present we are not able
to extend this result to larger times. 
Nevertheless, an interesting 
Lyapunov function is inherited by invoking the analogy of our
problem at infinite time with the PDE based Optimal Transport setting
and its relationship with the shape optimization problem.
Assuming existence of the solution and its first time derivative at
all times, we prove that the Lyapunov-candidate function is always
decreasing in time and reaches stationarity at $t=\infty$. Although we
are not currently able to use it properly, this result
seems promising in the search for a global existence theorem and a
formal justification of the conjecture that the proposed slime-mold
problem is, at infinite time, equivalent to the MK problem.

\subsection{Notations}

We will denote by $\mathcal{F}$ the set of essentially bounded
functions $f$ with zero mean and compact support in an open, bounded,
and convex and simply-connected subset $\Omega$ of $\REAL^n$: 
\begin{equation*}
  \mathcal{F}:=\left\{ 
    f\in L^{\infty}(\Omega): 
    \Supp(f)\subsetneq\Omega \mbox{ and } \int_{\Omega }f\dx=0
    \right\}.
\end{equation*}
Without loss of generality, we may redefine the domain of $f$ by means
of an $n$-dimensional ball $B(0,\bar{R})$ with $\bar{R}$ sufficiently
large so that  
$\Omega\subset B(0,\bar{R})$ and with $f=0$ outside of $\Omega$ 
(in most cases of interest $\Omega$ can be an $n$-dimensional interval).
We will still use the same symbol $\Omega$ to denote such subset.
We define the subset $\mathcal{D}$ of $\Cac$ as:
\begin{equation*}
  \mathcal{D}:=\left\{ \tdens\in\Cac \mbox{ such that } 
    \lambda(\tdens):=\min_{x\in\bar{\Omega}}\tdens(x)\ge\alpha>0\right\},
\end{equation*}
where
\begin{gather*}
  \Cac=\left\{v:\bar{\Omega}\mapsto\REAL\ :\ v_{[\delta,\bar{\Omega}]}:=
    \underset{x\neq y}{\sup} \frac{|v(x)-v(y)|}{|x-y|^{\delta}}<+\infty\right\}\\
  \|v\|_{\Cac}:=\underset{\bar{\Omega}}{\sup}\  v+v_{[\delta,\bar{\Omega}]}
\end{gather*}
with $\delta \in ]0,1[$.
We indicate with $\Hac$ the H\"older space of continuously
differentiable functions with first derivatives in $\Cac$.
In the proof of Theorem~\ref{exist:ODE}, we will be using the following
characterization of the subspace $\mathcal{D}$.
Given $0<a<b<\infty$, the subspace $\mathcal{D}$ can be re-defined as
the union of open and convex subsets:  
\begin{gather}
\label{union:d(a,b)}
 \mathcal{D}=\underset{0<a<b<+\infty}{\bigcup} \mathcal{D}(a,b)\\
 \nonumber
 \mathcal{D}(a,b):=\left\{ \tdens\in\Cac \mbox{ such that } 
    a<\lambda(\tdens)\le\|\tdens\|_{\Cac}<b\right\}.
\end{gather}

Standard results on regularity theory of elliptic
PDEs~\citep{giaquinta} allow us to give the following definitions.
\begin{Defin}[Potential]
  \label{def:U}
  Let $\tdens\in\mathcal{D}$ and $f\in\mathcal{F}$.
  The operator $\Pot:\mathcal{D}\mapsto \Hac$ 
  maps $\tdens$ into $\Pot(\tdens)$,
  where $\Pot(\tdens)$ is the unique weak solution of the elliptic
  equation~\eqref{sys:intro:div}, i.e.
  \begin{align*}
    &\int_\Omega \tdens\,\nabla u\cdot\nabla\varphi\dx 
      = \int_\Omega f\varphi\dx \qquad 
      \forall\varphi\in H^1(\Omega),\\    
    & \int_\Omega u\dx = 0.
  \end{align*}
\end{Defin}
\begin{Defin}[Flux]
  \label{def:B}
  Let $\tdens\in\mathcal{D}$ and $f\in\mathcal{F}$.
  The operator
  $\Flux:\mathcal{D}\mapsto \Cac$ is defined as:
  \begin{equation*}
    \tdens\mapsto \Flux(\tdens):=\tdens|\nabla \Pot(\tdens)|.
  \end{equation*}
\end{Defin}

\subsection{Local existence and a Lyapunov-candidate function}

In this section we report the main results and developments that lead
to the local well-posedness of the model and to the identification of
the proposed Lyapunov-candidate function.  
We start by writing the weak formulation of system~\eqref{sys:intro}:
\begin{subequations}
  \label{sys:weak}
  \begin{gather}
    \int_{\Omega} \tdens(t,x)\nabla  u(t,x)\nabla \varphi(x)\dx=
         \int_{\Omega}f(x)\varphi(x)\dx \quad 
          \forall \varphi \in H^1(\Omega) \label{sys:div:weak}\\ 
    \int_{\Omega} u(t,x)\dx=0 \label{bou:cond:u:weak}\\
    \tdens'(t,x)= \tdens(t,x)|\nabla u(t,x)|-\tdens(t,x)\label{sys:dyn:weak}\\
    \tdens(0,x)=\tdens_0(x)\in \mathcal{D}\label{bound:cond:d:weak} 	
  \end{gather}
\end{subequations}
Definitions~\ref{def:U} and~\ref{def:B} allow us to substitute the
potential $\Pot(\tdens)$ into~\eqref{sys:dyn:weak} to obtain the following
semilinear evolution equation: 
\begin{gather}
  \label{CP}
    \begin{array}{cc}
      \tdens'(t)=-\tdens(t)+\Flux(\tdens(t)) \\
      \tdens(0)=\tdens_0\in \mathcal{D}
    \end{array}
\end{gather}
where the dependence on $x$ has been omitted.  Obviously, the pair
$(\tdens(t,x), u(t,x))$ can be reconstructed from
$(\tdens(t),\Pot(\tdens(t)))$.  The mild formulation of~\eqref{CP} is
given by:
\begin{equation}
  \label{mild}
  \tdens(t)=e^{-t}\tdens_0+\int_0^{t}e^{s-t}\Flux\;(\tdens(s))\ds
\end{equation}
The Banach-Caccioppoli fixed point theory states that local existence
and uniqueness of~\eqref{CP} require that the operator
$\Flux(\tdens)$ be locally Lipschitz continuous. 
We have the following:
\begin{Prop}
  \label{pot_flux_lip}
  The Potential and Flux operators $\Pot$ and $\Flux$ are Lipschitz
  continuous and bounded in $\mathcal{D}(a,b)$ for all $a$ and $b$,
  $0<a<b<\infty$. 
\end{Prop}
In other words we have that  for every $\tdens\in\mathcal{D}(a,b)$,
\begin{align}
  \label{locbound}
  \|\Pot(\tdens)\|_{\Hac}&\le C_1(a,b) \\
  \|\Flux(\tdens)\|_{\Cac}&\le C_2(a,b)
\end{align}
and there exist constants $L_{\Pot}(a,b)$ and $L_{\Flux}(a,b)$ such
that, for every $\tdens_1,\tdens_2\in\mathcal{D}(a,b)$:
\begin{align}
  \label{loclip}
  \|\Pot(\tdens_1)-\Pot(\tdens_2)\|_{\Hac}
  &
    \le  L_{\Pot}(a,b)\|\tdens_1-\tdens_2\|_{\Cac} \\
  \|\Flux(\tdens_1)-\Flux(\tdens_2)\|_{\Cac}
  &
    \le L_{\Flux}(a,b)\|\tdens_1-\tdens_2\|_{\Cac}
  \end{align}
The previous Proposition follows from the following:
\begin{Lemma}
  \label{holder:f:G}
  Let $\tdens\in \mathcal{D}$, $F_0\in\mathcal{F}$, and
  $F=(F_i)_{(i=1,\ldots,n)}\in [\Cac]^{n}$ and
  $u\in H^1(\Omega)$ be the unique solution of 
  \begin{equation}
    \label{div:form}
    \left\{
      \begin{array}{cc}
        \int_{\Omega}\tdens\nabla u\cdot\nabla\varphi\dx 
        =\int_{\Omega}\left(F_0\varphi+F\cdot\nabla\varphi\right)\dx 
        \quad \forall \varphi\in H^1(\Omega)\\
        \int_{\Omega}u\dx =0
      \end{array}
    \right.
  \end{equation}
  Then $u\in \mathcal{C}^{1,\delta}(\bar{\Omega})$ and the following
  estimate holds:
  \begin{equation}
    \label{grad:u:Cac}
    \|\nabla u\|_{\Cac}\leq \Kdelta(n,\Omega,\delta)
    \Kmu(\tdens)\left(\|F_0\|_{L^{\infty}(\Omega)}+\|F\|_{\Cac}\right)
  \end{equation}
  where $\Kdelta(n,\Omega,\delta)$ is a constant depending on the 
  dimension $n$, the domain $\Omega$, and the H\"older regularity
  $\delta$ of $\tdens$, and:
  \begin{equation}
    \label{estimate:C(d)}
    \Kmu(\tdens)=\Kmu\left(\lambda(\tdens),\|\tdens\|_{\Cac}\right)=
    \frac{1}{\lambda(\tdens)}
    \left(
      \frac{\|\tdens\|_{\Cac}}{\lambda(\tdens)}
    \right)^{\frac{n+\delta}{2\delta}}.
  \end{equation}  
\end{Lemma}
This lemma, whose proof is given in Section~\ref{proof:holder:f:G},
extends classical results of regularity theory of elliptic equations
with H\"older continuous coefficients by careful estimation of the
dependence upon $\tdens\in\mathcal{D}$ of the constants $\Kdelta$ and
$\Kmu$~\citep{troi,giaquinta}.  The latter result allows us to prove
the main theorem of the paper:
\begin{Theo}
  \label{exist:ODE}
  Given $\tdens_0\in\mathcal{D}$ there exists $\tau(\tdens_0)>0$ such
  that the Cauchy Problem~\eqref{CP} admits a unique solution
  $\tdens\in \mathcal{C}^1\left([0,\tau(\tdens_0)[\,;\Cac\right)$.
  Moreover, this solution remains strictly greater than zero and:
  \begin{equation}
  \label{always:pos}
    \lambda(\tdens(t))\geq e^{-t}\lambda(\tdens_0)
  \end{equation}
  for $t\in[0,\tau(\tdens_0)[$.
\end{Theo}
This theorem suggests that there must be an interplay between $\tdens$
and $\nabla\Pot$ that constrains the flux $|\Flux|$ to remain bounded
so that, under the hypotheses of Theorem~\eqref{exist:ODE}, existence
and uniqueness of the solution pair $(\tdens,u)$ is expected for all
times. However, such a global result seems out of reach within the
current framework because the conclusions of Lemma~\ref{holder:f:G} do
not allow us to improve the local Lipschitzianity of $\Pot$ and
$\Flux$.  Nonetheless, local existence and uniqueness of the solution,
albeit incomplete, offer a certain degree of confidence on the
correctness of the model of the PP dynamics, and justify the use of
numerical discretizations to supply credible evidence that the model
is well-posed for all times.

A further indication that the problem is globally well-posed is
provided by the fact that we have identified a Lyapunov-candidate
function that could be used to constrain the problem.  This function
is derived by borrowing from results in the field of
mass/shape optimization and considering its relationship with optimal
transport problems.  We start the derivation of the Lyapunov-candidate
function by defining the shape optimization problem as described
in~\citet{shape}.

Assume to have two nonnegative functions $f^+$ and $f^-$ in $\Omega$,
with $\int_{\Omega}f^+=\int_{\Omega}f^-$, representing, e.g., the
density of a given electric charge, and a fixed amount of a conductive
material, described by a nonnegative function $\sigma$ having unit
mass.  Interpreting $\Omega$ as an insulating medium, we ask the
question how to optimally distribute the given conductive material so
that the heating induced by the current flow from $f^+$ to $f^-$ is
minimal. In~\citet{shape} the authors consider the case where
$f^+$, $f^-$ and $\sigma$ are non-negative Radon measures and discuss the
connection between the mass optimization problem and a generalized
version of the MK equations.  They proved that the optimal
distribution of conductive material $\sigma^*$ given $f^+$ and $ f^-$
is the normalized optimal transport density $\mu^*$ of the MK equations
with $f=f^+-f^-$
In the case of $f\in \mathcal{F}$, these results together with those
given in~\citet{DePascale:2004} lead to the following formulation: 
\begin{displaymath}
  \min_{\sigma \in \mathcal{M}} \mathcal{E}(\sigma)
     = \frac{1}{2}\int_{\Omega}\sigma|\nabla \Pot (\sigma)|^2
  \qquad \mathcal{M}=\left\{\sigma \in L^{\infty}(\Omega) \ : 
    \int_{\Omega}\sigma \dx =1\right\}
\end{displaymath}
Recasting the shape optimization
intuition into the framework of our formulation, it is natural to
identify a density
$\sigma(t)=\left({\tdens(t)}/{\int_{\Omega}\tdens(t)\dx}\right)$ that
belongs to $\mathcal{M}$.  Then, following our conjecture that
$\tdens(t)\rightarrow\tdens^*$ as $t\rightarrow\infty$, it is natural
to assume that $\sigma(t)$ should tend to $\sigma^*$, solution of the
shape optimization problem.  
Noting that
$\Pot(\sigma)=\left(\int_{\Omega}\tdens\dx\right)\Pot(\tdens)$, 
we can define our Lyapunov-candidate function as:
\begin{equation}
  \label{def:L}
  \Lyap(\tdens):= \frac{1}{2}\int_{\Omega}\tdens\dx \cdot  
     \int_{\Omega}\tdens|\nabla \Pot(\tdens)|^2\dx
\end{equation}
Intuitively, we are looking for a density $\tdens$
that gives the best trade off between the total mass, and thus the
cost, of the transport infrastructure and the dissipated energy. 
In the large times limit, the above function $\Lyap$ evaluated along the
trajectory $\tdens(t)$ should tend to the minimum of~\eqref{def:L} and
its time derivative should always be negative.
Thus, $\Lyap$ is a reasonable Lyapunov-candidate function, as the
following theorem states: 
\begin{Theo}
\label{theo:lyap}
  The function $\Lyap:\mathcal{D}\;\longmapsto\REAL^+$ defined above
  is strictly decreasing in time along the solution $\tdens(t)$ for
  $t\in[0,\tau(\tdens_0)[$. 
\end{Theo}
The previous theorem leads to the following lemma:
\begin{Lemma}
 \label{L1:bound}
 For all $t\in [0,\tau(\tdens_0)[$ the 
 $L^1$-norm of $\tdens(t)$ and $\Flux(\tdens(t))$ are bounded
 by constants depending only on the initial data
 $\tdens_0$.
\end{Lemma}

\section{Numerical solution}

In this section we describe the numerical discretization used to find
an approximate solution to the proposed model, and report the numerical
results obtained solving the proposed model to simulate the Physarum
Polycephalum dynamics on the maze.
Next we report some numerical results aimed at justifying the
conjecture that the solution of the proposed time-dependent model
tends as time tends to infinity to the solution of the
Monge-Kantorovich problem. To this aim, the model is applied to the
solution of the OT problems proposed by~\citet{prigozhin} and the
numerical results are compared.

\subsection{Spatial and temporal discretization}

The numerical solution of eqs.~\eqref{sys:weak} is obtained by means
of the method of lines.  We denote with $\Triang[h](\Omega)$ and
$\Triang[2h](\Omega)$ two regular triangulations of the domain
$\Omega$, where $\Triang[h]$ is obtained from $\Triang[2h]$ by
connecting the edge mid-points. Indicating with $N_{2h}$, $E_{2h}$,
and $M_{2h}$ the number of nodes, edges, and triangles of mesh
$\Triang[2h]$, respectively, we have that $N_{h}=N_{2h}+E_{2h}$ and
$M_{h}=4M_{2h}$.  Spatial discretization is achieved using standard
Galerkin finite elements by projecting eq.~\eqref{sys:div:weak} onto
the finite dimensional space
$V_h=\PONE(\Triang[h](\Omega))=\Span\{\varphi_{1}(x),\ldots,\varphi_{N_h}(x)\}$
of piecewise linear Lagrangian basis functions defined on
$\Triang[h](\Omega)$ and eq.~\eqref{sys:dyn:weak} onto the finite
dimensional space $W_{h}=\PZERO(\Triang[2h](\Omega))$ of piecewise
constant functions.  
Following this approach, the spatially discrete
potential $u_h(t,x)$ and diffusion coefficient $\tdens_h(t,x)$ are
written as:
\begin{gather*}
  u_{h}(t,x)=\sum_{l=1}^{N_{h}} u_{l}(t)\varphi_{l}(x) \qquad
  \varphi_{l}\in V_h=\PONE(\Triang[h]) \\
  \tdens_h (t,x)= \sum_{r=1}^{M_{2h}}\tdens_{r}(t) \chi_{r}(x) \qquad
  \chi_{r}(x)=\left\{
    \begin{array}{ll}
      1 \ \text{if} \ x\in T_{r}\\
      0 \ \text{if} \ x\notin T_{r}
    \end{array}
    T_{r}\in\Triang[2h] \right.
\end{gather*}
Thus the FEM method yields the following system of differential
algebraic equations:

Find $(u_h,\tdens_h)\in V_h\times W_h$ such that:
\begin{subequations}
  \label{fem}
  \begin{gather}
    a_{\tdens_h}(u_{h},\varphi_{m})=
      \int _{\Omega}\tdens_{h}\nabla u_{h}\cdot\nabla\varphi_{m}\dx
    =
      \int_{\Omega}f\varphi_{m}\dx\qquad m=1,\ldots,N_{h}
      \label{fem:elliptic}\\
    \int_{\Omega} u_h \dx = 0 \label{fem:zero:mean} \\
    \int _{\Omega}\tdens'_{h}\chi_{s}\dx
    = 
      \int_{\Omega}\tdens_{h}\left(|\nabla u_{h}|-k\right)\chi_{s} \dx
      \qquad s=1,\ldots,M_{2h}
      \label{fem:ode} \\
    \int_{\Omega}\tdens_{h}(0,\cdot)\chi_{s}\dx=
       \int_{\Omega}\tdens_{0}\chi_{s}\dx
      \qquad s=1,\ldots,M_{2h}
    \end{gather}
\end{subequations}
Forward Euler time stepping discretizes eq.~\eqref{fem:ode}
yielding a decoupled system of linear equations. 
Denoting with $\Deltat$ the time-step size so
that $t_{j+1}=t_{j}+\Deltat$ and writing $u_{h}^{j}\approx u_{h}(t_j,x)$ and
$\tdens_{h}^{j}\approx \tdens_{h}(t_j,x)$, 
the final solution algorithm can be written as:
\begin{subequations}
  \label{algorithm}
  \begin{gather}
    a_{\tdens_{h}^{j}}(u_h^{j},\varphi_m)=(f,\varphi_m), 
        \qquad m=1,\ldots,N_{h} \label{alg:div}\\
    \tdens_{s}^{j+1}= 
    \tdens_{s}^{j}\left(1+\Deltat(|\nabla u_h^{j}|_{s}-k_{s})\right),
        \qquad s=1,\ldots,M_{2h}  \label{alg:dyn}
  \end{gather}
\end{subequations}
where 
\begin{displaymath}
    a_{\tdens_{h}^{j}}(u_{h}^{j},\varphi_{m})=
      \int_{\Omega}\tdens_{h}^{j}\nabla
      u_{h}^{j}\cdot\nabla\varphi_{m}\dx =
      \sum_{l=1}^{N_{h}}u_{l}^{j}
      \int_{\Omega}\tdens_{h}^{j}\nabla\varphi_{l}\cdot\nabla\varphi_{m}\dx
\end{displaymath}
and
\begin{displaymath}
  |\nabla u_{h}|_{s} = \frac{1}{|T_{s}|}\int_{T_{s}}|\nabla u_{h}|\dx 
  \qquad
  k_{s} = \frac{1}{|T_{s}|}\int_{T_{s}}k\dx 
  \qquad
  \tdens_{s}^0=\frac{1}{|T_{s}|}\int_{T_{s}}\tdens_{0}\dx.
\end{displaymath}
Eq.~\eqref{alg:div} is a sparse system of linear equations of
dimension $N_{h}$, and is solved by means of the incomplete Choleski
preconditioned conjugate gradient with the variant proposed
by~\citet{Bochev:2005} to solve the corresponding semi-definite linear
system arising in our pure Neumann problem.  The
$M_{2h}\times M_{2h}$ linear system described in eq.~\eqref{fem:ode}
is diagonal and leads to~\eqref{alg:dyn}.
To maintain the coercivity of the FEM bilinear form
$a_{\tdens_{h}^{j}}(\cdot,\cdot)$ we impose a lower bound on
$\tdens_h$ of $10^{-10}$. 

\begin{remark}
  The choice of the two different FEM spaces $V_{h}$ and $W_{h}$,
  which resemble the inf-sup stable $\PONE$-iso-$\PTWO$/$\PONE$  FEM
  spaces for the Stokes equation~\citep{brezzi}, is dictated by the
  need to reduce or eliminate oscillations on the gradient $|\nabla
  u_h|$ that occur if different approximating spaces are used. 
  Experimentally, in fact, we observed that in the case
  $(u_{h},\tdens_{h})\in\PONE(\Triang[h])\times\PZERO(\Triang[h])$
  oscillations in $\nabla u_{h}$ appear and destroy convergence of the
  numerical solution. 
  On the other hand, using the proposed approach
  $(u_{h},\tdens_{h})\in\PONE(\Triang[h])\times\PZERO(\Triang[2h])$,
  which essentially is nothing else than a simple average of the
  solution gradient on the larger triangles, all the oscillations
  disappear. 
\end{remark}

\subsection{Numerical simulation of Physarum Polycephalum dynamics}

\begin{figure}
  \centerline{
    \includegraphics[width=0.45\textwidth]{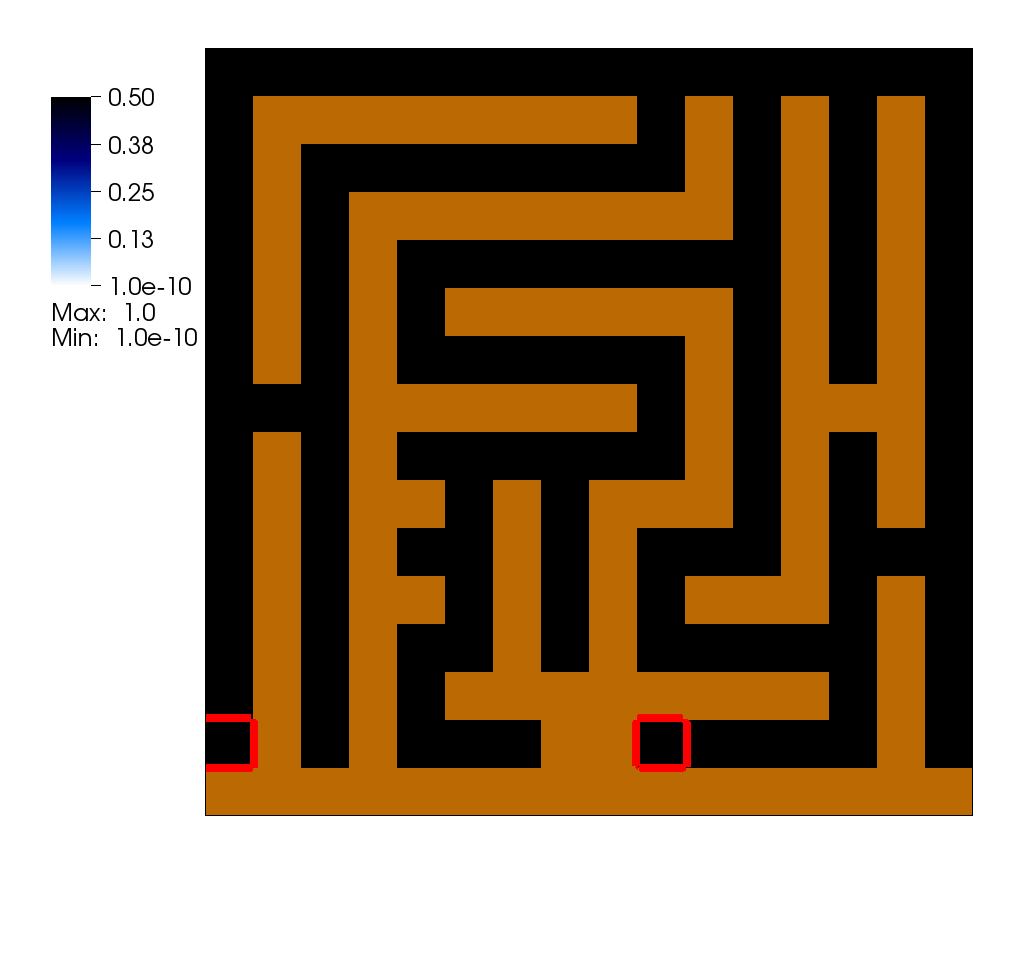}
    \includegraphics[width=0.45\textwidth]{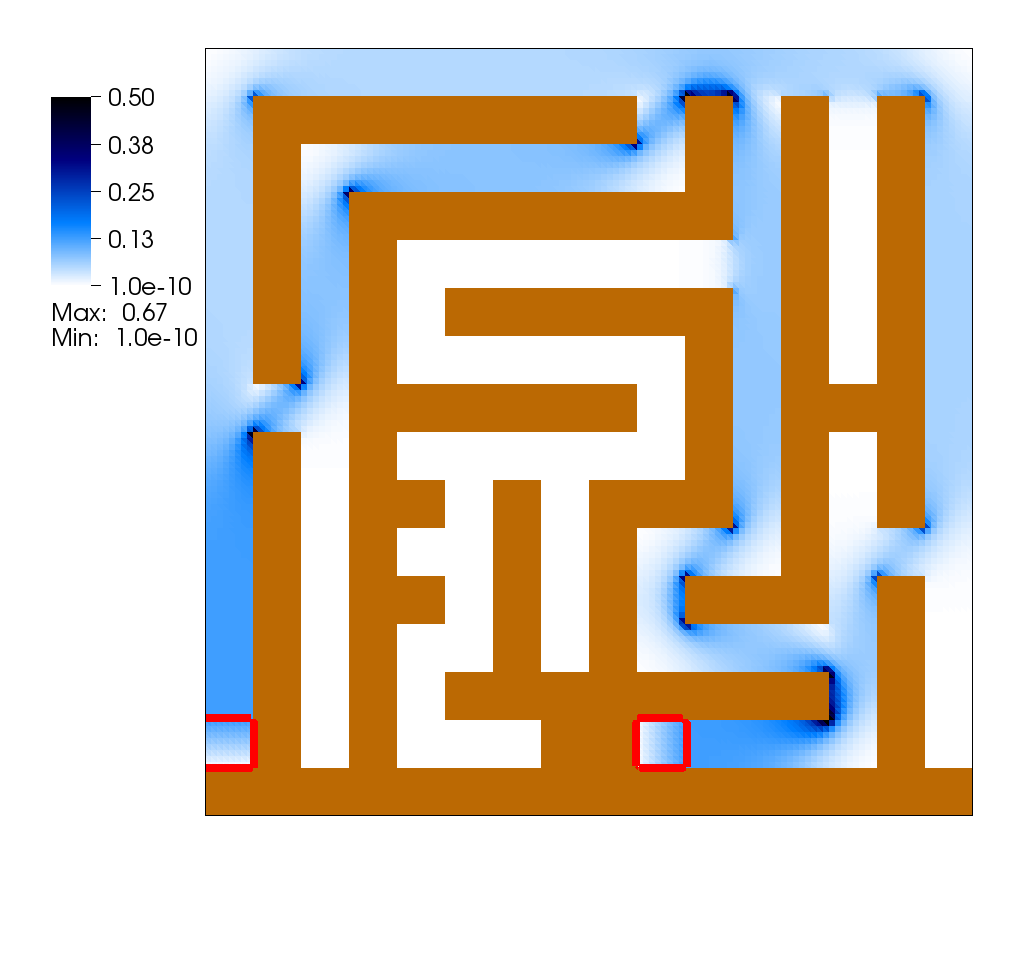}
  }
  \centerline{
    \includegraphics[width=0.45\textwidth]{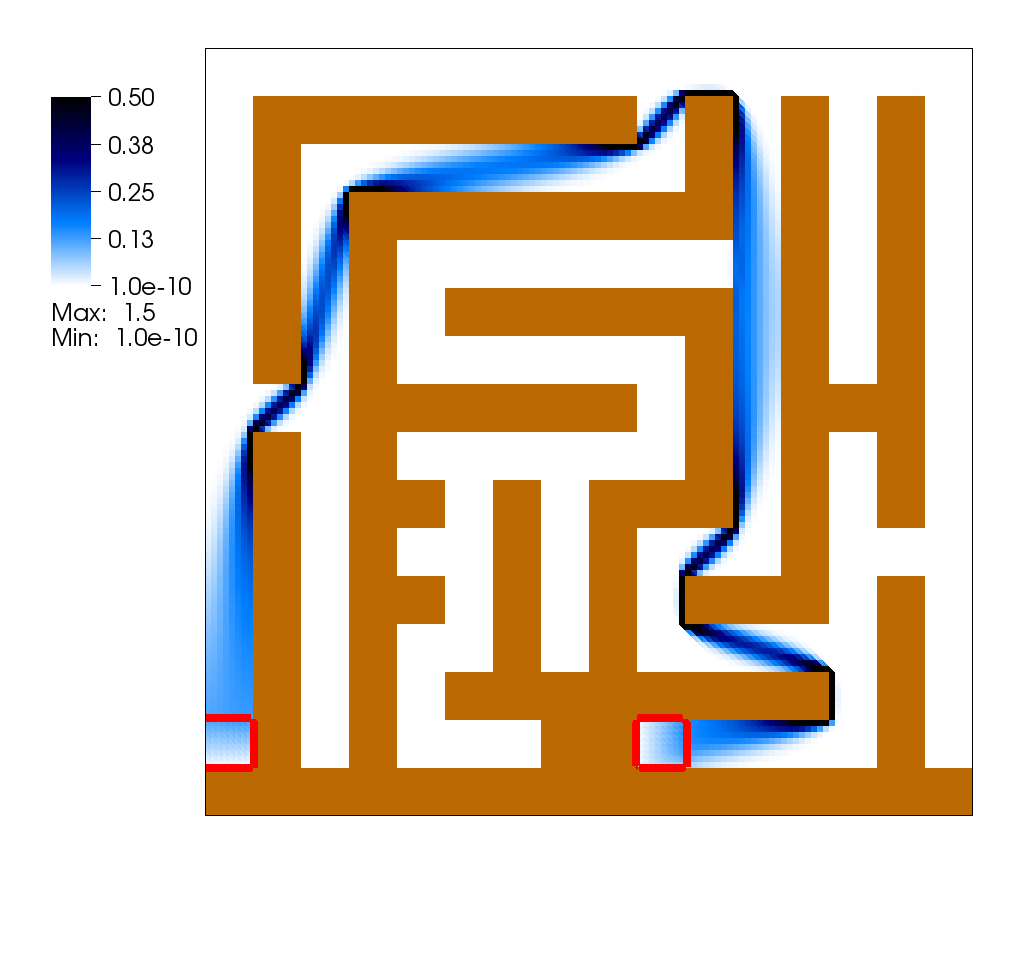}
    \includegraphics[width=0.45\textwidth]{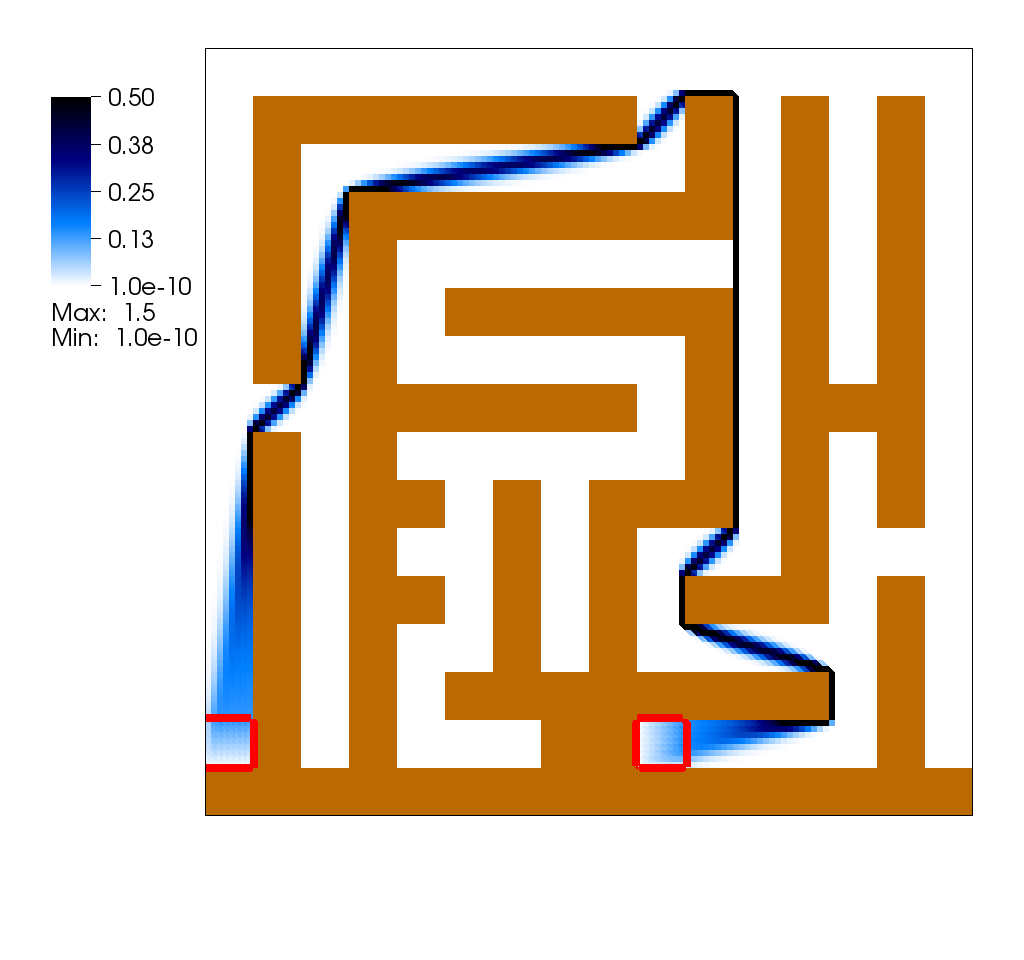}
  }
  \caption{Simulation of the dynamics of Physarum Polycephalum mass
    reorganization in the maze experiment of~\citet{nakagaki:maze}:
    distribution of PP transport density at 
    dimensionless times $t=0$ (top, left), $t=5.96$ (top, right),
    $t=60.3$ (bottom, left), and $t=9.6\times 10^3$ (bottom, right). The
    simulation was done on a triangulation $\Triang[2h]$ with 32768
    triangles and 16641 nodes. At the last time step of the simulation
    the $\tdens_h$ variation was smaller than $\tau=5\times 10^{-9}$.}
  \label{fig:maze-dynamics}
\end{figure}

The proposed model and its numerical discretization just described are
applied to the simulation of the dynamics of the experiment described
in~\citet{nakagaki:maze}. The domain encompassing the entire maze
setup is discretized by means of a uniform triangulation obtained by
subdividing each edge of the square-shaped maze into 128
subdivisions yielding the coarser mesh $\Triang[2h]$ comprised of 32768
elements and 16641 nodes.  This high resolution is required to follow
accurately the walls of the maze, which are described by setting
$k(x)=1000$ (brown colors in the upper left panel of
Figure~\ref{fig:maze-dynamics}) while the maze paths are characterized
by $k(x)=1$.  The initial condition $\tdens_0$ is set to $10^{-10}$ on
the maze walls and one elsewhere.  The two food sources $f^+=1$ and
$f^-=-1$ are shown as red squares in the figure.  We employ a variable
time step size starting from $\Delta t_0=10^{-2}$ and with 
$\Delta t_{j+1}=\min(1.01\Deltat,0.5)$, to ensure stability of the
Forward Euler scheme is verified for all times with an ample safety margin.
The simulation is stopped when the relative difference in
$\tdens_h$ becomes smaller than the tolerance $\tau$, i.e.:
\begin{equation}\label{variation}
    \var(\tdens_h(t^{j})) = \frac{\|\tdens_h^{j+1}-\tdens_h^{j}\|_{L^2(\Omega)}}
    {\Deltat\|\tdens_h^{j}\|_{L^2(\Omega)}}\le \tau
\end{equation}
with $\tau=5\times{10}^{-9}$.
Figure~\ref{fig:maze-dynamics} shows the distribution of
$\tdens_{h}$ at three
different times, in addition to the initial data. The three times are
chosen in agreement with the simulations reported
in~\citet{tero:model} and are used to highlight the intermediate
phases when the P.Polycephalum starts retreating from the dead ends
($t=5.96$) and when it starts to concentrate along the identified
shortest path ($t=60.3$). The final steady state configuration is
achieved at $t=9.6\times 10^{3}$. 
Note that the same numerical solution is obtained, albeit at different 
dimensionless times, starting from different initial conditions
$\tdens_0$. 
Looking at the first intermediate phase we see that PP has retreated
from the dead end paths of the maze but persists on all the possible
paths connecting the two food sources. We note a stronger
concentration of $\tdens_{h}$ at the edges of the maze walls, indicating
that PP starts accumulating around a narrow band along the shortest
route. This is clear from the solution at the next time (lower left
panel), when the shortest path within the maze is completely
identified but $\tdens_{h}$ is still somewhat dispersed. At the final
time, $\tdens_{h}$ is distributed along the optimal route displaying
varying approximation levels depending on the alignment of the mesh
triangles with the support of $\tdens_{h}$. In fact, the vertical
portion is one element thick, while the oblique routes encompass more
than one triangle. 
All these observations are in line with the results proposed
by~\citet{tero:model}, although in our case the graph structure is not
imposed a-priori but it is mimicked through the appropriate definition
of $k(x)$. 
Note that in our continuous setting the presence of $k(x)$ is related
to the cost of through-flow, while the original graph-based
formulation does allows only flow throuh the graph edges. 

\subsection{Numerical simulation of optimal transport problems}

\begin{figure}
  \centerline{
    \includegraphics[width=0.4\textwidth]{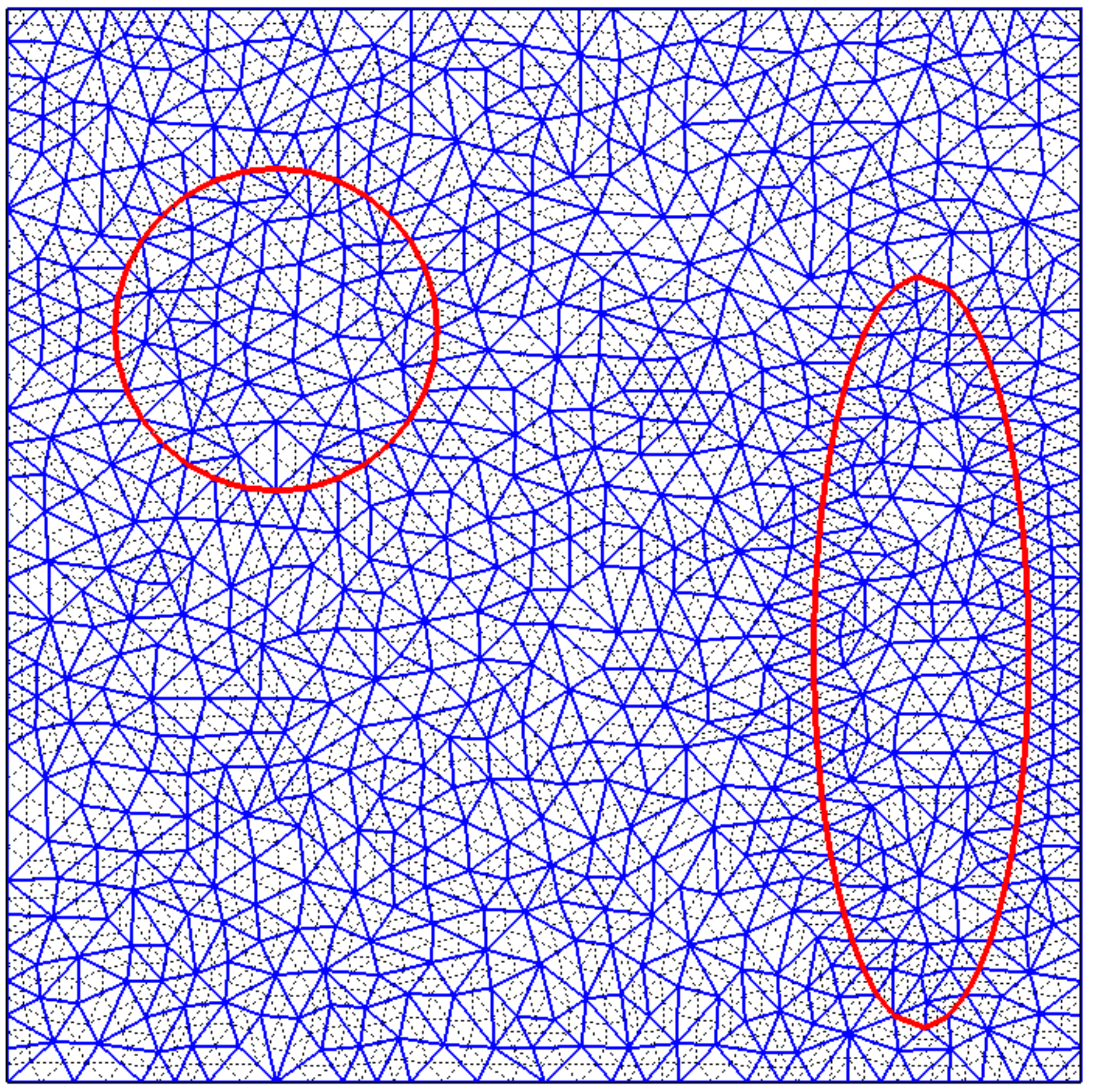}
    \includegraphics[width=0.4\textwidth]{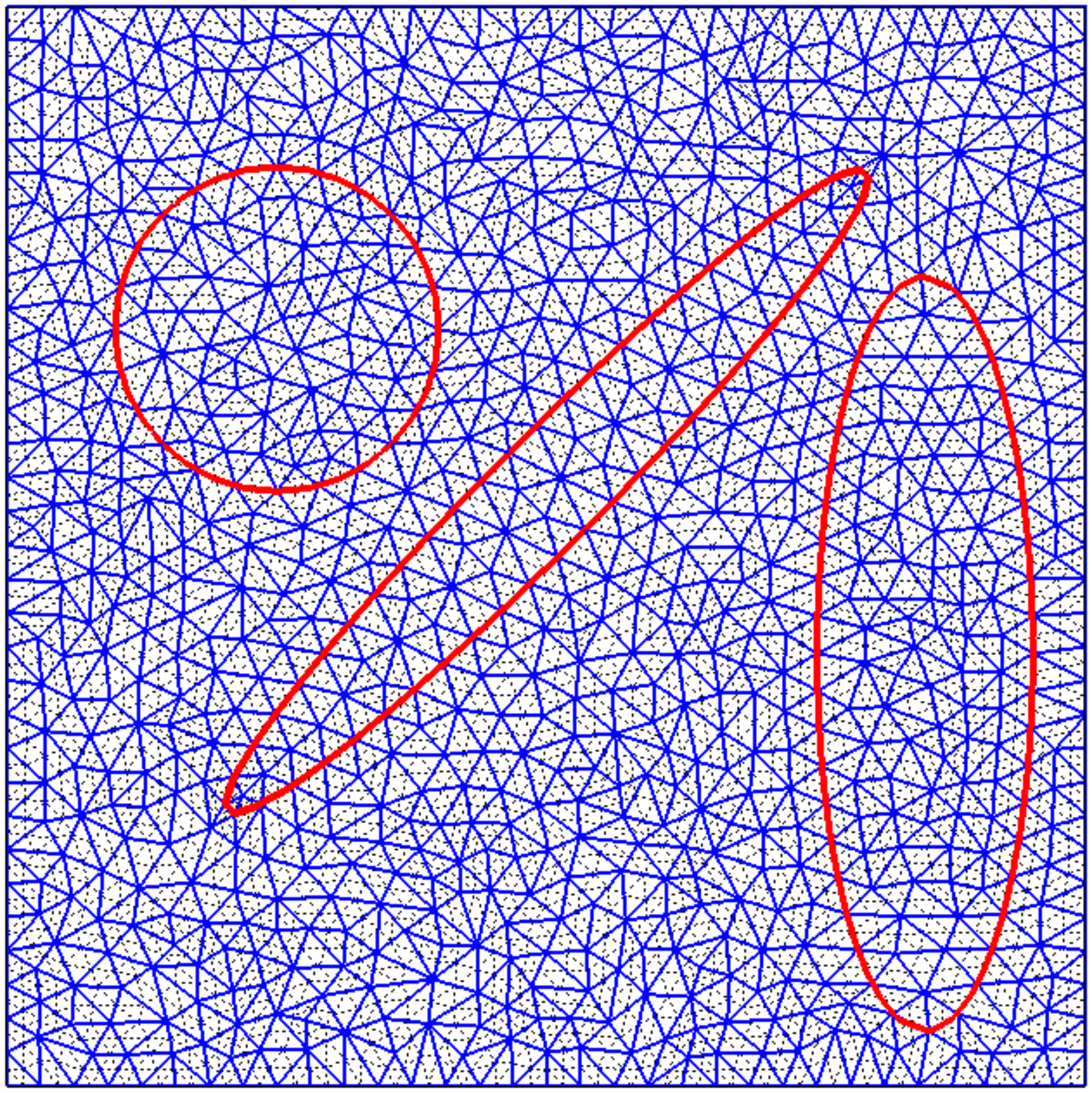}
  }
  \caption{Unit square domain and triangulations $\Triang[2h]$ (solid
    blue) and $\Triang[h]$ (dotted black) for the homogeneous (left)
    and heterogeneous (right) test cases. The supports of 
    $f^+$ (left circle) and $f^-$ (right ellipse) are shown in red,
    as well as the central ellipse where $k(x)$ is different from the 
    background in the heterogeneous case.}
  \label{domain-prig}
\end{figure}

In this section we report numerical evidence supporting the conjecture
that the long-time solution of the slime-mold dynamic model is
solution of the Monge-Kantorovich equations~\eqref{evans}.
To this aim we compare our solution with the one proposed
by~\citet{prigozhin}. 
At the same time, the aim of this section is to show that the solution
of the MK equations by means of the proposed dynamic model is
numerically easier than the direct solution of system~\eqref{evans},
suggesting that the introduction of time relaxes the numerical
difficulties yielding a more robust, flexible, and efficient method
for the numerical solution of optimal transportation problems.

We apply the numerical approach described above to the solution of the
homogeneous and heterogeneous model problems proposed in Example 1
of~\citet{prigozhin}.  These problems consider a unit square domain
with zero Neumann boundary and constant forcing having uniform
positive value ($f^+$) supported on a circle and uniform negative
value ($f^-$) on an ellipse in such a way that spatial balance is
ensured (Figure~\ref{domain-prig}).  The domain contains a central
oblique oval shape where a value for $k(x)$ different from the unitary
background value is specified.  Four different test cases are
defined. The first one, called the homogeneous test case, considers a
unit value of $k(x)$ assigned in the entire domain. This corresponds
to the standard MK equation reported in system~\eqref{evans}. Then
three heterogeneous problems are addressed by setting the value of
$k(x)$ in the central oblique ellipse equal to $k_e=0.01$ and $k_e=100$,
respectively.  The former case favors flow through the central oval,
while the latter impedes it.  The final test case considers an
intermediate value of $k_e=3$ in the central oval and is used to show
the dynamical behavior of the proposed model by looking at the
numerical solution at intermediate times.

To accurately impose the forcings of the different problems, the
triangulations are adapted to the supports of $f(x)$ and $k(x)$,
compelling the use of nonuniform meshes.  The homogeneous case is
solved on a sequence of three uniformly refined triangulation couples
$\Triang[2h](\Triang[h])$. The coarsest mesh has $820(3170)$ nodes and
$1531(6124)$ triangles (Figure~\ref{domain-prig}, left), for a total
of $4701(=3170+1531)$ degrees of freedom in the final algebraic
solution algorithm~\eqref{algorithm}.  The next two levels have
$3170(12463)$ and $6124(24496)$ triangles, yielding a total of $18587$
degrees of freedom, and $12463(49421)$ nodes and $24496(97984)$
triangles, for a total of $73917$ degrees of freedom, respectively.
In the heterogeneous case we employ two triangulation levels starting
with a mesh of $933(3603)$ nodes and $1738(6952)$
elements~(Figure~\ref{domain-prig}, right).  All simulations start
with a spatially uniform unitary initial condition $\tdens_{h}^{0}$
and are run until condition~\eqref{variation} is satisfied with
$\tau=5\times{10}^{-9}$.
The sequence of time step sizes is the same as in the case of the
slime mold simulation, but with different upper bounds to ensure
stability of Euler scheme at large times.

\begin{figure}
  \centerline{
    \includegraphics[width=0.32\textwidth]{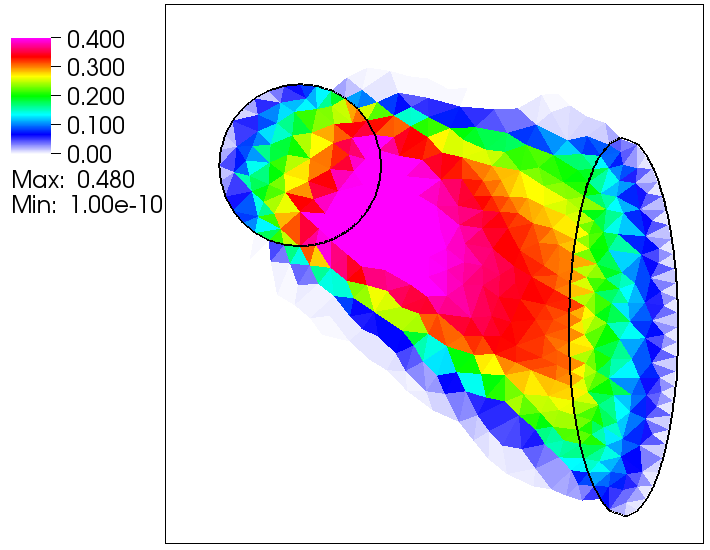}
    \includegraphics[width=0.32\textwidth]{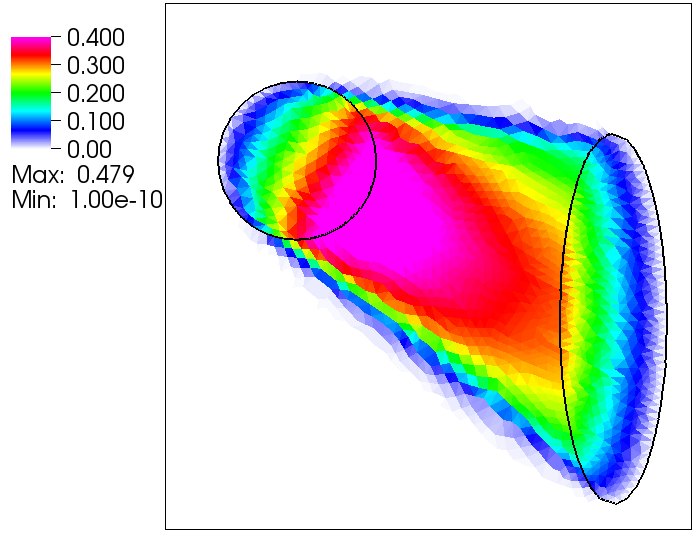}
    \includegraphics[width=0.32\textwidth]{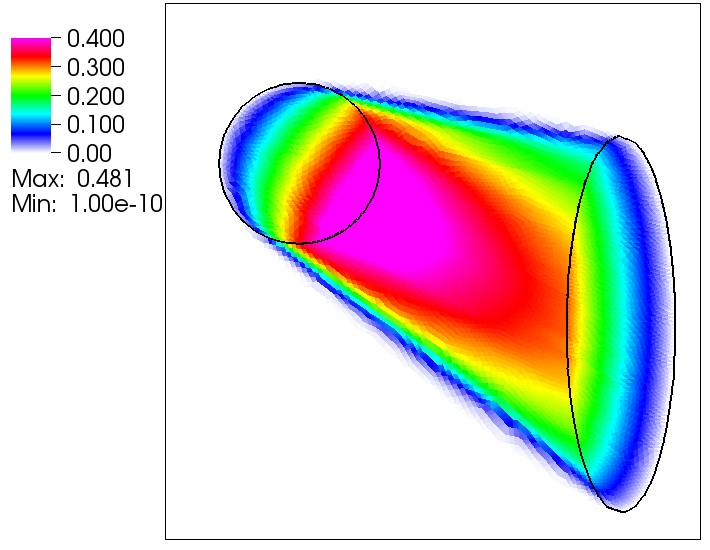}
  }
  \caption{ Numerical solution of the homogeneous OT
    problem: distribution of the transport density
    $\tdens_h$ at three different mesh refinement levels of the
    coarser ($\Triang[2h]$) triangulation (from left to right).}
  \label{num_res}
\end{figure}

\begin{figure}
  \centerline{
    \includegraphics[width=0.5\textwidth]{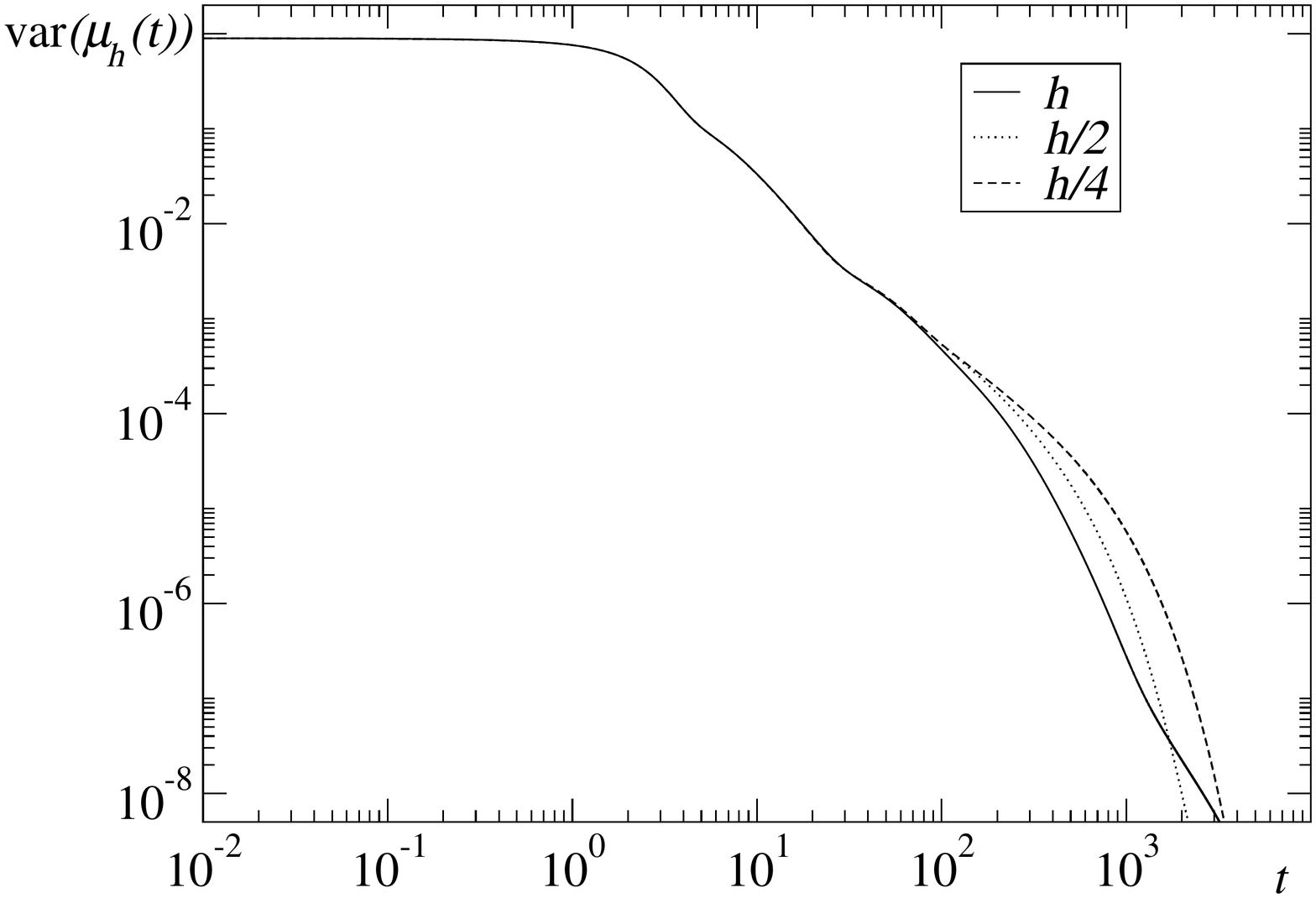}\qquad
    \includegraphics[width=0.5\textwidth]{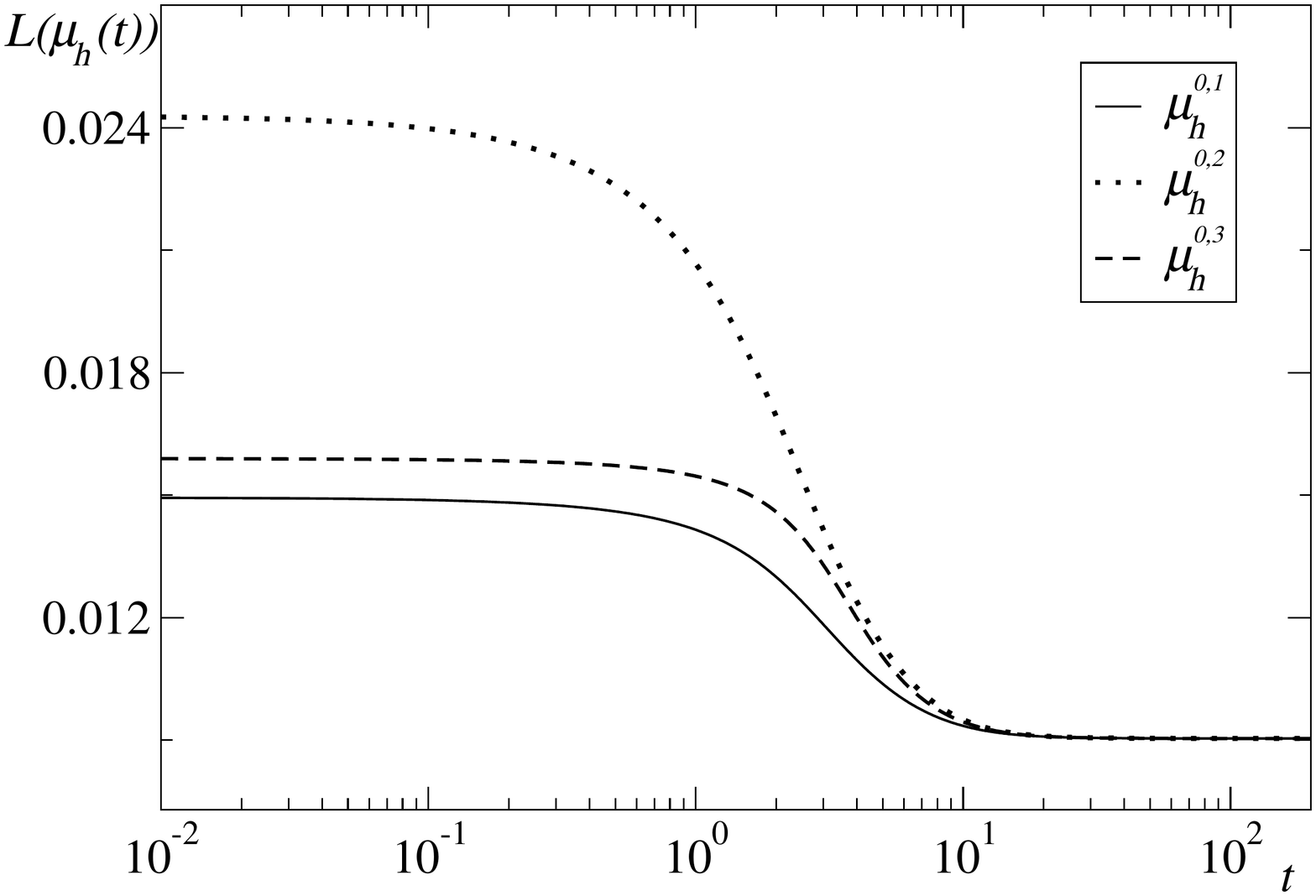} 
  }
  \caption{Numerical results for the homogeneous OT problem.  Left
    panel: convergence towards zero of the relative variation of the
    transport density
    $\var(\tdens_h(t^{j})) = \|\tdens_h^{j+1}-\tdens_h^{j}\|_{L^2(\Omega)}/
    (\Deltat\|\tdens_h^{j}\|_{L^2(\Omega)})$ for the three refinement
    levels.  
    Right panel: behavior of the Lyapunov function
    $\Lyap(\tdens_{h}(t^{j}))$ vs time 
    for the three different initial data ($\tdens_{h,1}^{0}$,
    $\tdens_{h,2}^{0}$, $\tdens_{h,3}^{0}$) and the coarsest mesh level.}  
  \label{num_res_lyap}
\end{figure}

The experimental results for the homogeneous test case are shown in
Figure~\ref{num_res}, where the spatial distributions of
$\tdens_{h}^{j}$ and $|\nabla u_{h}^{j}|$ are plotted for the coarser
($\Triang[h]$) mesh at the three refinement levels.  All the features
of the expected solution of the MK problem are present in these
results, which are in good agreement with those of~\citet{prigozhin}.
The support of the transport density concentrates on the region
connecting the circle boundaries with the ellipse.  Along the
transport rays, $\tdens_{h}$ increases from zero to its maximum value
within the supports of $f^+$ and then decreases moving towards the
ellipse, tending to negligible values at ray
ends~\citep{evans2,buttazzo}.  The norm of the gradient, not shown for
brevity, is practically one in all triangles of $\Triang[2h]$ lying
within the support of $\tdens_h$.  All these features are clearly
visible already at the coarsest level, although the effects of large
elements is evident. At increasing mesh refinements the delineation of
the support of $\tdens_{h}$ is sharper, and at the final level the
accuracy seem to be satisfactory and compares well with the solution
proposed by~\citet{prigozhin}.

Convergence toward steady state is monotonic, as shown in
Figure~\ref{num_res_lyap} (left), where the relative variation of
$\var(\tdens_{h}(t))$ (eq.~\eqref{variation}) is shown in log-log
scale as a function of time for the three refinement levels.  We
notice that the convergence behavior is identical for all refinement
levels up to a relative $\tdens_{h}$ variation of approximately
$10^{-3}$, which corresponds to a dimensionless time
$\hat{t}\approx 60$. At this point the discretization errors
seem to slightly influence convergence. We note that from this time on
the value of the Lyapunov candidate function
(Figure~\ref{num_res_lyap}, right panel) is practically constant,
suggesting that the solution has effectively converged to a stationary
state.
At the first refinement level, the number of time steps to reach
$\hat{t}$ is 431, while 6641 is the number of time steps at the final
time $t^*=3200$ (the upper bound on $\Deltat$ is 0.25).  
The average number of iterations of the IC0-preconditioned conjugate
gradient per time-step is 66. 

The model solution at large times is insensitive to the distribution
of $\tdens_{h}^{0}$ as shown by the behavior of the Lyapunov function
reported in the right panel of Figure~\ref{num_res_lyap}, where three
different sets initial data, $\tdens_{h}^{0,i}$, $i=1,2,3$, are
tested:
\begin{align*}
  \tdens_{h,1}^{0}(x,y)=1
  \quad
  \tdens_{h,2}^{0}(x,y)=
  0.1+
  4\left(
    (x-0.5)^2+(y-0.5)^2
  \right) \\
  \tdens_{h,3}^{0}(x,y)=
  3+
  2 \sin(8 \pi x)\sin(8\pi y)
\end{align*}
This result is upheld by the fact that the Lyapunov function
$\Lyap(\tdens_h)$ numerically evaluated starting from three different
initial conditions always converges to the same value
(Fig.~\ref{num_res_lyap}, right).  Its non-increasing behavior shows
that $\Lyap$ is a plausible candidate for a Lyapunov function.

In this work we are not interested in implementing the most
computational efficient algorithm but we want to show that simple
numerical methods are sufficient to find an accurate solution to the
MK equations. We would like to remark that many improvements can be
done to the numerical scheme and are indeed under development.
Notwithstanding its simplicity, our approach is competitive with
respect to more direct MK solution methods, such as the one proposed
in~\citet{prigozhin}.  These authors solve the direct optimal transport
problem using a mixed finite element method with adaptive mesh
refinement that converged to a final triangular grid adapted to the
shape of the transport plan. The discretization on the final mesh
level leads to a final nonlinear system with approximately 60000
degrees of freedom which is solved by an ad hoc nonlinear successive
over-relaxation method.  The Successive Over-relaxation method
used to solve the nonlinear system was considered converged when the
relative flux residual was smaller than $10^{-3}$.
Noting that in this case the convergence criteria based on relative
variations of $\tdens_h$ or on $q=\tdens_h|\nabla u_h|$ are
equivalent, as stated above, from Figure~\ref{num_res_lyap} we see
that this convergence level is reached at $\hat{t}\approx 60$ in our
case. At this time, our candidate Lyapunov function has reached an
almost steady condition with a very small rate of decrease, signaling
that for practical purposes convergence to the sought solution has been
achieved. 

\begin{figure}
  \centerline{
    \includegraphics[width=0.45\textwidth]{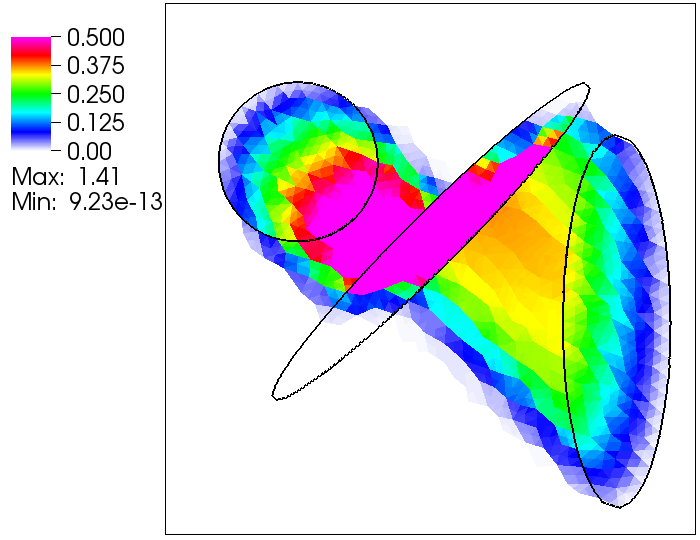}
    \includegraphics[width=0.45\textwidth]{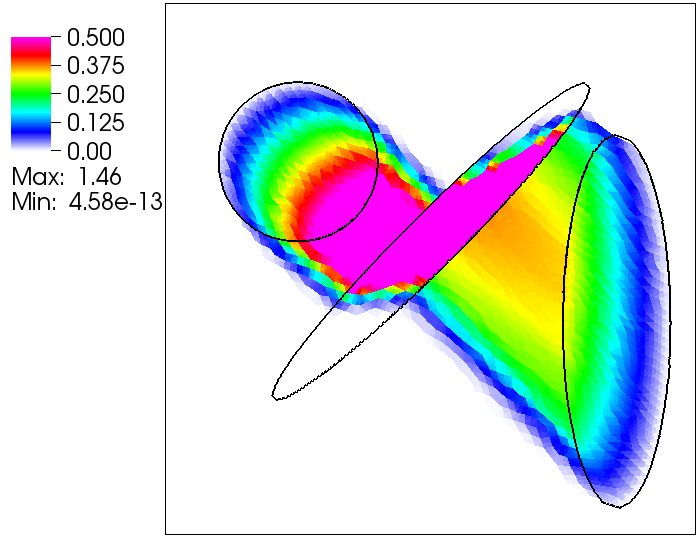}
  }
  \centerline{
    \includegraphics[width=0.45\textwidth]{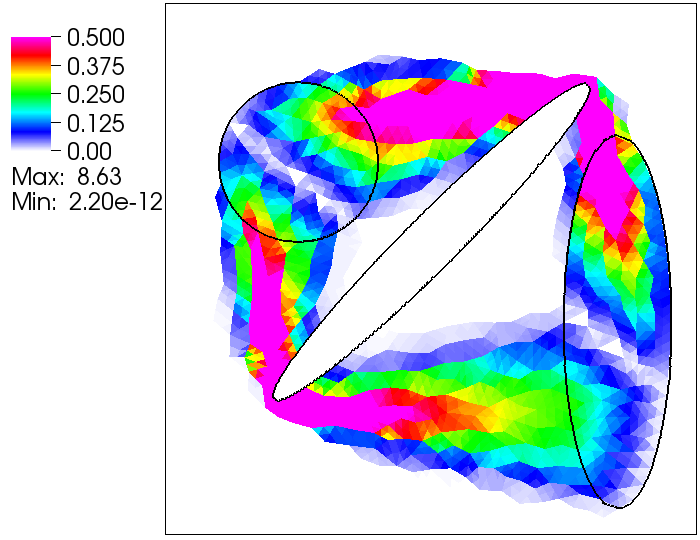}
    \includegraphics[width=0.45\textwidth]{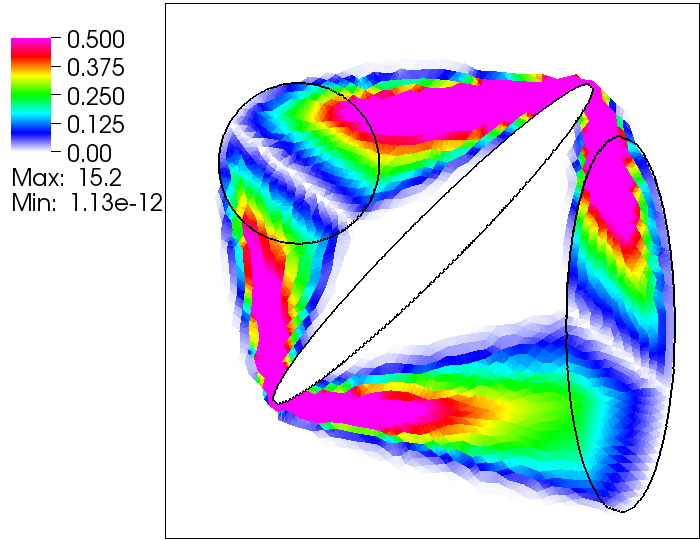}
  }
  \caption{ Numerical solution of the heterogeneous test case for
    $k_e=0.01$ (top) and $k_e=100$ (bottom) in the central
    ellipse. The figures show the optimal transport flux
    $\tdens_h|\nabla u_h|$ for the mesh with 1738(6952) triangles and
    933(3603) (left) and the once-refined mesh 6952(27808) triangles
    and 3603(14157) nodes (right).}
  \label{num_res_hetero_k1}
\end{figure}

The numerical results for the heterogeneous case are shown in
Figure~\ref{num_res_hetero_k1} where the steady-state spatial
distribution of the absolute flux $q=\tdens_h|\nabla u|$ is plotted
in the case of $k_e=0.01$ (top panels) and $k_e=100$ (bottom panels). 
We first note that in this heterogeneous case the gradient is bounded
by $k(x)$ and not by one as in the previous test case. For this
reason we chose to to plot the flux $|q|=\tdens_h|\nabla u_h|$ instead
of $\tdens_h$.  
Two successively refined triangulations are used, leading to linear
systems of dimensions $N_{h}+M_{2h}=3603+1738$ and
$N_{h}+M_{2h}=14157+6952+1738$ for the coarser and the finer meshes,
respectively.  
The results are qualitatively comparable with those
of~\citet{prigozhin}, although obtained with a much coarser
discretization. 
It is interesting to note that the qualitative features of the
solution are obtained already at the coarser mesh, with no visible
numerical artifacts barring mesh roughness. We would like to stress
here the fact that, notwithstanding the fact that the mesh nodes are
not aligned with the support of $\tdens_h$, the geometrical features
of the solution are well captured at all mesh resolution levels. 
From the flux spatial distribution, we see that values of $k_e$ lower
than one promotes larger fluxes across the central ellipses. On the
contrary, values substantially larger than one restricts through-flow, and
promotes the circumnavigation of the high low conductivity areas.

\begin{figure}
  \centerline{
    \includegraphics[width=0.32\textwidth]{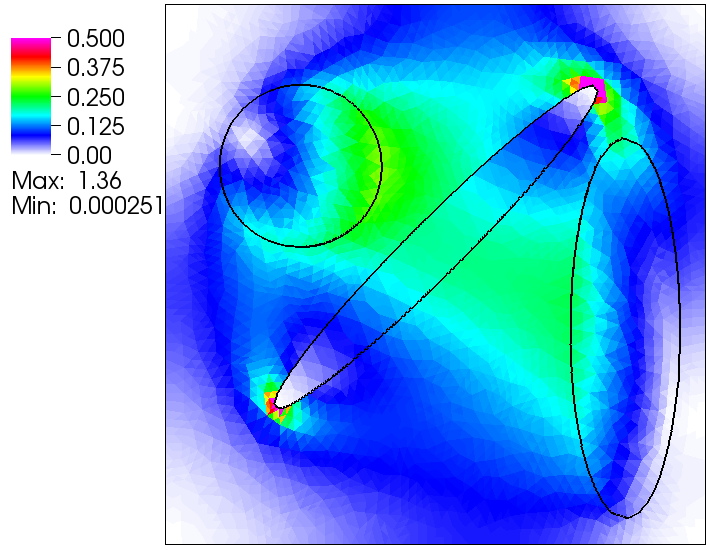}
    \includegraphics[width=0.32\textwidth]{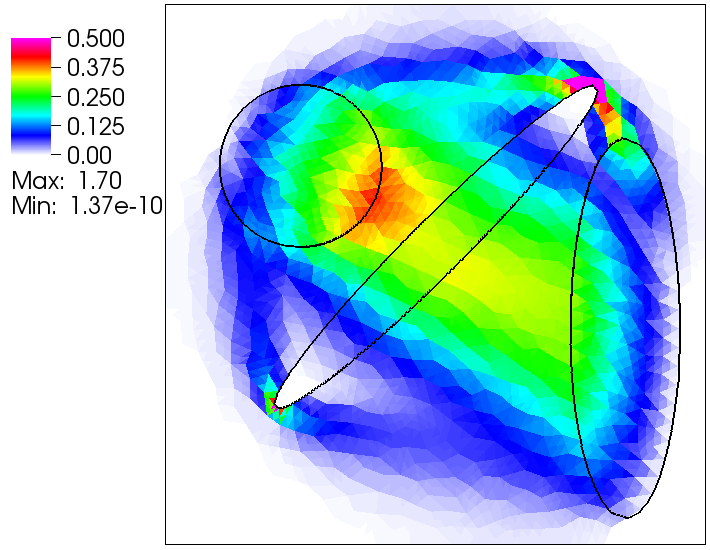}
    \includegraphics[width=0.32\textwidth]{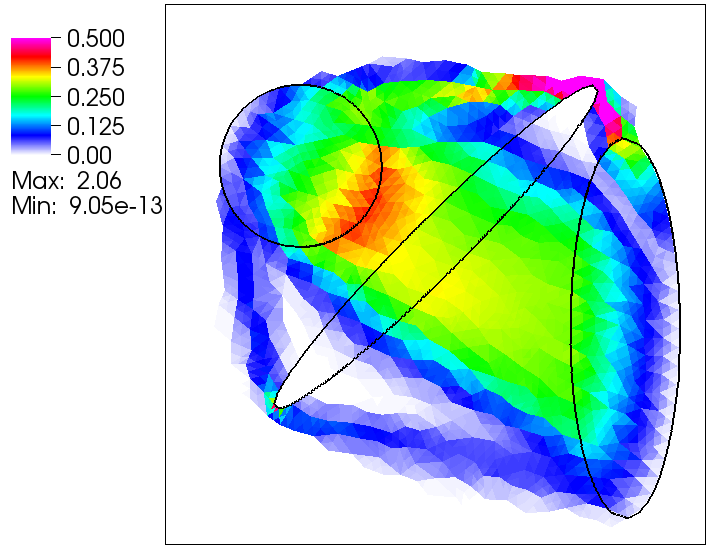}
  }
  \centerline{
    \includegraphics[width=0.32\textwidth]{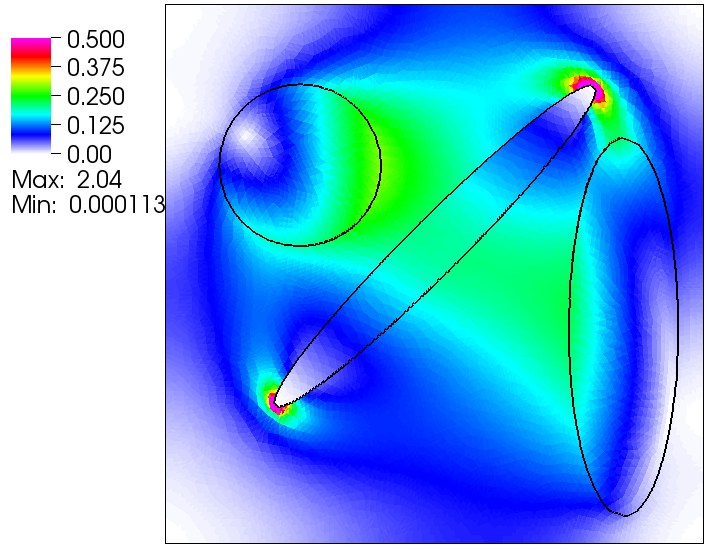}
    \includegraphics[width=0.32\textwidth]{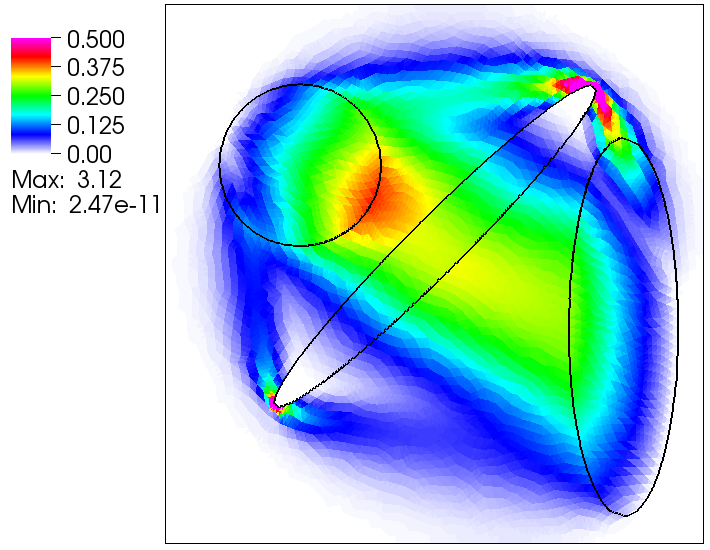}
    \includegraphics[width=0.32\textwidth]{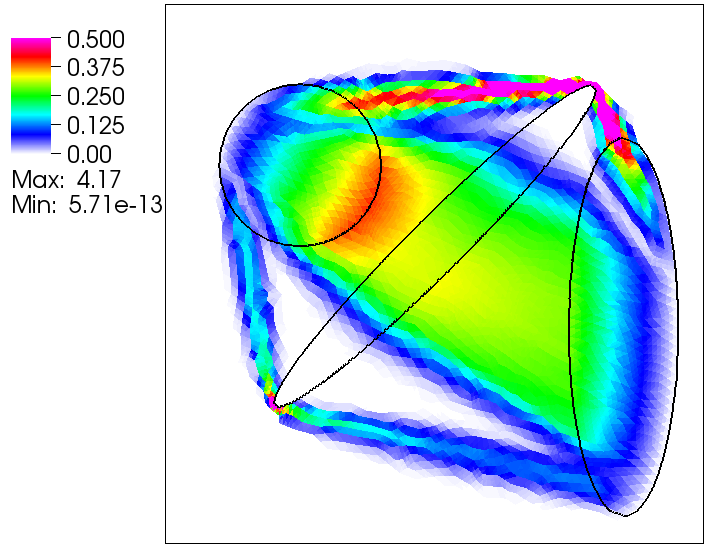}
  }
  \caption{
    Numerical solution of the heterogeneous test case with $k_e=3$
    in the central ellipse in terms of 
    optimal transport flux $\tdens_h|\nabla u_h|$ 
    for 3 different times ($\var(\tdens_h)=0.1, 0.01, 5\times 10^{-9}$
    from left to right) and 2 refinement levels (from top
    to bottom).
  }
  \label{num_res_hetero_k3}
\end{figure}

Finally, Figure~\ref{num_res_hetero_k3} shows the distribution of
$\tdens_h$ for the case of $k_{e}=3$ at three different times (left to
right) and two successive refinement levels (top to bottom).  The
three times are chosen so that the $\tdens_{h}$ variation reaches the
thresholds $\var(\tdens_{h}(\hat{t}_1))=0.1$,
$\var(\tdens_{h}(\hat{t}_2))=0.01$, and
$\var(\tdens_{h}(\hat{t}_3))=5\times 10^{-9}$.  Correspondingly, we
have $\hat{t}_1\approx 5.2$ and $\hat{t}_2\approx 21$, remaining the
same for both tested triangulations, and $\hat{t}_3=1600$ for the
coarser level and $\hat{t}_3=2200$ for the finer mesh as steady state
is achieved at a later time for the finer mesh, reflecting the fact
that the overall error is smaller. In fact, the converged steady state
solution occurs after 6616 and 8955 time steps for the coarse and fine
triangulations, respectively.  Note that, for this last heterogeneous
test case, the time-stepping sequence employed an upper bound on
$\Deltat$ equal to 0.25.
Also in this case the steady-state numerical solution is similar to
the results reported by~\citet{prigozhin}. We see from the time
sequence that our model constructs the transport map gradually.
Starting from the uniformly distributed initial condition, it first
identifies the the larger flow paths and then refines them to arrive
at the final configuration. 
The overall results are consistently pointing towards the veridicity
of the conjecture that the infinite-time solution to our problem
indeed coincides with the solution of the Monge-Kantorovich equations
in the support of the optimal transport path.

\section{Proofs of the results}
\label{proofs}
\subsection{Proof of Proposition~\ref{pot_flux_lip}: Lipschitz
  continuity of $\Pot$ and $\Flux$}
\label{proof:pot_flux_lip}

We start by considering the solution $u\in H^1(\Omega)$ of
eq.~\eqref{sys:intro:div} in weak form: 
\begin{equation}
  \label{weak}
  \left\{
    \begin{array}{cc}
      a_\tdens(u,\varphi)=\int_{\Omega}\tdens\nabla u\cdot\nabla\varphi\dx
        =\int_{\Omega}f\varphi \quad \forall \varphi\in H^1(\Omega)\dx\\[1.5ex]
      \int_{\Omega}u\dx =0
    \end{array}
  \right.
\end{equation}
We recall here the general hypotheses for the existence and uniqueness
of the solution $u$ of~\eqref{weak}: i) the domain $\Omega$ 
be a bounded, connected, and convex subset of $\REAL^n$ with
$\mathcal{C}^1$ or Lipschitz 
boundary $\partial\Omega$; ii) $f\in L^2(\Omega)$; iii) 
the bilinear form $a_\tdens(u,\varphi)$ be bounded and coercive, i.e.:
\begin{gather}
  \label{bound}
  \exists\  0<\Lambda<+\infty \quad \mbox{such that} \quad 
  \arrowvert a_\tdens(u,v)\arrowvert \leq \Lambda \|u\|_{H^1(\Omega)}\|v\|_{H^1(\Omega)}
  \quad \forall u,v \in H^1(\Omega) \\[1.0em]
  \label{coe}
  \exists\  0<\lambda<+\infty \quad \mbox{such that} \quad  
      a_\tdens(u,u) \geq \lambda \|u\|^2_{H^1(\Omega) }
  \quad \forall u\in H^1(\Omega)
\end{gather}
It follows from the above that, if $\Omega =B(0,\bar{R})$, $\tdens \in
\mathcal{D}$, and $f\in\mathcal{F}$, the hypotheses for existence and
uniqueness of the elliptic equation are satisfied and thus the
operators $\Pot$ and $\Flux$ are well defined. 
The regularity result provided in Lemma~\ref{holder:f:G} ensures that
$\Pot$ maps into $\Hac$. 
Then, given $\tdens\in\mathcal{D}(a,b)$  
we apply Lemma~\ref{holder:f:G} with $F_0=f$ and $F_i=0$ to obtain:
\begin{equation}
  \label{Pot_bound}
  \|\nabla \Pot(\tdens)\|_{\Cac}\leq 
     \Kdelta(\Omega,n,\delta) \Kmu(a,b)\|f\|_{L^{\infty}(\Omega)}
\end{equation}
From this, the boundedness of the Potential $\Pot(\tdens)$ follows
immediately. 
The local Lipschitz continuity of $\Pot$ derives from the following
considerations. 
First, we note that given $\tdens_1,\tdens_2\in\mathcal{D}(a,b)$ and 
$u_k=\Pot(\tdens_k)$ with $k=1,2$, we have:
\begin{gather}
  \notag
  \int_{\Omega}\tdens_1\nabla u_1\nabla\varphi\dx=\int_{\Omega}f \varphi\dx =
  \int_{\Omega}\tdens_2\nabla u_2\nabla\varphi\dx \quad 
      \forall \varphi\in H^1(\Omega)\\
  \label{eq:lip}
  \int_{\Omega}\tdens_1\nabla(u_1-u_2)\nabla\varphi\dx=
  \int_{\Omega}(\tdens_2-\tdens_1)\nabla u_2\nabla\varphi\dx \quad 
      \forall \varphi \in H^1(\Omega)
\end{gather}
Then, application of Lemma~\ref{holder:f:G} with $F_0=0$ and
$F=-(\tdens_1-\tdens_2)\nabla\Pot(\tdens_2))$, which belongs to
$(\Cac)^n$ because of~\eqref{Pot_bound}, yields:
\begin{align}
  \notag \|\nabla(u_1-u_2)\|_{\Cac}\leq& 
     \Kdelta(\Omega,n,\delta)\Kmu(\tdens_1)
           \|(\tdens_1-\tdens_2)\nabla u_2\|_{\Cac} \\
  \notag \leq&
     \Kdelta(\Omega,n,\delta)\Kmu(\tdens_1)\|\tdens_1-\tdens_2\|_{\Cac}
        \|\nabla u_2\|_{\Cac}\\
  =&\Kdelta(\Omega,n,\delta)^2 \Kmu(a,b)^2
     \|f\|_{L^{\infty}(\Omega)}\|\tdens_1-\tdens_2\|_{\Cac}
     \label{Pot_lip}
\end{align}
We can also prove that the flux $\Flux$ is bounded in
$\mathcal{D}(a,b)$. 
In fact, since the H\"older norm is sub-multiplicative, we can write
from~\eqref{Pot_bound}: 
\begin{equation*}
  \|\Flux(\tdens)\|_{\Cac}=\|\tdens |\nabla\, \Pot(\tdens)|\;\|_{\Cac}\leq
     \Kdelta(\Omega,n,\delta)\; b\,\Kmu(a,b)\|f\|_{L^{\infty}(\Omega)}
\end{equation*}
Lipschitz continuity of $\Flux$ derives from~\eqref{Pot_lip} as follows:
\begin{align*}
  \|\Flux(\tdens_1)-\Flux(\tdens_2)\|_{\Cac}
  &=\|\tdens_1|\nabla\Pot(\tdens_1)|-
    \tdens_2|\nabla\Pot(\tdens_2)|\;\|_{\Cac}\\
  &=\|\tdens_1(|\nabla\Pot(\tdens_1)|-|\nabla\Pot(\tdens_2)|)
    -(\tdens_2-\tdens_1)|\nabla \Pot(\tdens_2)|\;\|_{\Cac}\\
  &\leq\|\tdens_1\|_{\Cac}\|\; 
    |\nabla [\Pot(\tdens_1)-\Pot(\tdens_2)]|\;\|_{\Cac} \\
  &\qquad  
    +\Kdelta(\Omega,n,\delta)\Kmu(a,b)\|f\|_{\infty}\|\tdens_1-\tdens_2\|_{\Cac}\\
  &\leq L_{\Flux}(a,b)\|\tdens_1-\tdens\|_{\Cac}
\end{align*}

\subsection{Proof of Theorem~\ref{exist:ODE}: local existence}
\label{proof:exist:ODE}

Given $\tdens_0\in\mathcal{D}$ and the Lipschitz continuity of $\Flux$
in $\mathcal{D}(a,b)$ asserted by Prop.~\ref{pot_flux_lip}, noting
that the subspace $\mathcal{D}$ can be decomposed as given
in~\eqref{union:d(a,b)}, we may restrict our investigations on 
$\mathcal{D}(a,b)$.
Standard arguments in the theory of ODEs in Banach Spaces
ensure local existence and uniqueness of the solution $\tdens(t)$.
In other words, there exists a sufficiently small $\tau(\tdens_0)>0$
such that the fix point problem~\eqref{mild}
admits a solution $\tdens\in\mathcal{C}^0([0,\tau(\tdens_0);\Cac)$.
The Lipschitz continuity of $\Flux$ automatically ensures that
$\tdens\in\mathcal{C}^1([0,\tau(\tdens_0);\Cac)$. 
The proof of estimate~\eqref{always:pos} is immediate by considering
that in~\eqref{mild} the term $\int_0^{t}e^{s-t}\Flux(\tdens(s))\ds$ 
is always greater than zero.



\subsection{Proof of Theorem~\ref{theo:lyap}: Lyapunov-candidate
  function} 
\label{proof:lyap}

In this section we report the proof that Lyapunov-candidate
function $\Lyap$ decreases along the trajectories.
Before we compute the time derivative along the trajectories 
we need to prove the $\mathcal{C}^1$-regularity in time of $u(t)$.
We have the following (proof in Section~\ref{sec:proof-regularity}):
\begin{Prop}
  \label{u:C1}
  The function $u(t)$ belongs to the space
  $\mathcal{C}^1([0,\tau(\tdens_0)[;\mathcal{C}^{1,\delta}(\Omega))$.
  For each $t\in[0,\tau(\tdens_0)[$ its time derivative $u'(t)$ solves
  the following equation: 
  \begin{equation}
    \label{u'}
    \left\{
      \begin{array}{cc}  
        \int_{\Omega}\tdens(t)\nabla u'(t)\cdot\nabla\varphi\dx =
       -\int_{\Omega}\tdens'(t)\nabla u(t)\cdot\nabla\varphi\dx \quad 
        \forall\varphi\in H^1(\Omega)\\
        \int_{\Omega}u'(t)\dx=0
      \end{array}
    \right.
  \end{equation}
\end{Prop}
We can now compute the time derivative of $\Lyap(\tdens(t))$ and prove it
is strictly negative, thus proving Theorem~\ref{theo:lyap}.
In fact we have:
\begin{align*}
 2 \Lyap'(t)
  &=\frac{d}{dt}\left(\int_{\Omega}\tdens(t)\dx
    \int_{\Omega}\tdens(t)|\nabla u(t)|^2\dx\right)\\
  &=\int_{\Omega}\tdens'(t)\dx\int_{\Omega}\tdens(t)|\nabla u(t)|^2\dx 
    +\int_{\Omega}\tdens(t)\dx\frac{d}{dt}\int_{\Omega}\tdens(t)|\nabla u(t)|^2\dx\\
  &=\int_{\Omega}\tdens'(t)\dx\int_{\Omega}\tdens(t)|\nabla u(t)|^2\dx
    +\int_{\Omega}\tdens(t)\dx\int_{\Omega}[\tdens'(t)|\nabla u(t)|^2 \\
  &\qquad\qquad  +2\;\tdens(t)\nabla u'(t)\cdot\nabla u(t)]\dx
\end{align*}
Substituting $\varphi=u(t)$ in \eqref{u'} we get:
\begin{displaymath}
 \int_{\Omega}\tdens(t)\nabla u'(t)\cdot \nabla u(t)\dx=
    -\int_{\Omega} \tdens'(t)|\nabla u(t)|^2\dx
\end{displaymath}
Thus
\begin{align*}
  2\Lyap'(t)
  &=\int_{\Omega}\tdens'(t)\dx\int_{\Omega}\tdens(t)|\nabla u(t)|^2\dx+
    \int_{\Omega}\tdens(t)\dx\int_{\Omega}[\tdens'(t)|\nabla u(t)|^2 \\
  &\qquad\qquad -2\;\tdens'(t)|\nabla u(t)|^2]\dx\\
  &=\int_{\Omega}\tdens'(t)\dx\int_{\Omega}\tdens(t)|\nabla u(t)|^2\dx
    -\int_{\Omega} \tdens(t)\int_{\Omega} \tdens'(t)|\nabla u(t)|^2\dx\\
  &=\int_{\Omega}\tdens(t)(|\nabla u(t)|-1)\dx
    \int_{\Omega}\tdens(t)|\nabla u(t)|^2\dx\\
  &\qquad\qquad
    -\int_{\Omega}\tdens(t)\dx
    \int_{\Omega}\tdens(t)(|\nabla u(t)|-1)|\nabla u(t)|^2\dx\\
  &=\int_{\Omega}\tdens(t)|\nabla u(t)|\dx
    \int_{\Omega}\tdens(t)|\nabla u(t)|^2\dx-
    \int_{\Omega}\tdens(t)\dx\int_{\Omega}\tdens(t)|\nabla u(t)|^3\dx
\end{align*}
Rewriting the product of integrals as an integral in
$\Omega\times\Omega$  of the functions
\begin{gather*}
  f(t;x,y):=\tdens(t,x)|\nabla u(t,x)|\;\tdens(t,y)|\nabla u(t,y)|^2-
    \tdens(t,x)\;\tdens(t,y)|\nabla u(t,y)|^3\\
  g(t;x,y):=\tdens(t,y)|\nabla u(t,y)|\;\tdens(t,x)|\nabla u(t,x)|^2-
    \tdens(t,y)\;\tdens(t,x)|\nabla u(t,x)|^3
\end{gather*}
we obtain:
\begin{align*}
  2\Lyap'(t)
  &=\int_{\Omega \times \Omega}f(t;x,y)\dx\, dy=
    \int_{\Omega \times \Omega}g(t;x,y)\dx\, dy\\
  &=\int_{\Omega \times \Omega}\frac{f(t;x,y)+g(t;x,y)}{2}\dx\, dy\\
  &=\int_{\Omega \times \Omega}\tdens(t,x)\;\tdens(t,y)\; 
    \frac{h(t;x,y)}{2}\dx dy
\end{align*}
where
\begin{align*}
 h(t;x,y)
 &=|\nabla u(t,x)|\ |\nabla u(t,y)|^2-|\nabla u(t,y)|^3
   +|\nabla u(t,x)|^2 |\nabla u(t,y)|\\
 &\qquad\qquad 
   -|\nabla u(t,x)|^3\\
 &=|\nabla u(t,x)|\; |\nabla u(t,y)|\left(|\nabla u(t,x)|
   + |\nabla u(t,y)|\right)\\
 &\qquad\qquad
   -\left(|\nabla u(t,x)|^3 +|\nabla u(t,y)|^3\right)\\
 &=|\nabla u(t,x)|\; |\nabla u(t,y)|\left(|\nabla u(t,x)|
   + |\nabla u(t,y)|\right)\\
 &\qquad\qquad
   -\left(|\nabla u(t,x)| +|\nabla u(t,y)|\right)\left(|\nabla u(t,x)|^2 
    +|\nabla u(t,y)|^2 \right.\\
 &\qquad\qquad   
   \left.-|\nabla u(t,x)|\; |\nabla u(t,y)|\right)\\
 &=-\left(|\nabla u(t,x)|+ |\nabla u(t,y)|\right)
    \left(|\nabla u(t,x)|^2 +|\nabla u(t,y)|^2\right.\\
 &\qquad\qquad
   \left.-2 |\nabla u(t,x)|\; |\nabla u(t,y)|\right)\\
 &=-\left(|\nabla u(t,x)|+ |\nabla u(t,y)|\right)\left(|\nabla u(t,x)|- 
    |\nabla u(t,y)|\right)^2
\end{align*}
Finally, we arrive at:
\begin{align*}
  \Lyap'(t)&=-\frac{1}{2}\int_{\Omega \times \Omega} \tdens(t,x)\tdens(t,y)
    \left(|\nabla u(t,x)|+|\nabla u(t,y)|\right)\left(|\nabla u(t,x)|- 
    |\nabla u(t,y)|\right)^2\dx\dy \\
  & \leq 0
\end{align*}

\begin{remark}
  Note that, assuming global existence of the solution, $\Lyap'=0$
  only if $\tdens=0$ or $|\nabla u|=\text{const}$ for all times. In
  particular, this assertion provides further support to our
  conjecture.
\end{remark}

\subsection{Proof of Lemma~\ref{L1:bound}: boundedness of $\tdens$ and
  $\Flux$}
\label{proof:L1:bound}

The $L^1$-bound for $\Flux(\tdens(t))$ derives directly from
Theorem~\eqref{theo:lyap}.
By Cauchy-Schwartz inequality, we can write:
\begin{gather*}
  \int_{\Omega}\Flux(\tdens(t)) \dx 
  = \int_{\Omega}\tdens(t)^{\frac{1}{2}}\,\tdens(t)^{\frac{1}{2}}|\nabla\Pot(\tdens(t))|\dx 
  \leq 
     \left(
       \int_{\Omega}\tdens(t)\dx\;\int_{\Omega}\tdens(t)|\nabla\Pot(\tdens(t))|^2\dx
     \right)^{\frac{1}{2}}\\
  \leq
     \left(
       \int_{\Omega}\tdens(0)\dx\;\int_{\Omega}\tdens(0)|\nabla\Pot(\tdens(0))|^2\dx 
     \right)^{\frac{1}{2}}
  =\Lyap(\tdens(0))^{\frac{1}{2}}
\end{gather*}
Then the $L^1$-bound for $\tdens(t)$ is obtained by integrating 
both sides of eq.~\eqref{mild} over $\Omega$:
\begin{align*}
  \int_{\Omega}\tdens(t)\dx
  &=\int_{\Omega}\left(e^{-t}\tdens_0+\int_0^{t}e^{s-t}\Flux\;(\tdens(s))\;ds\right)\dx\\
  &=e^{-t}\int_{\Omega}\tdens_0\dx+e^{-t}
      \int_0^{t}e^{s}\left(\int_{\Omega}\Flux\;(\tdens(s))\dx \right)\; ds\\
  &\leq e^{-t}\int_{\Omega}\tdens_0\dx+e^{-t}\int_0^{t} e^{s}(\Lyap(\tdens_0))^{\frac{1}{2}}\; ds\\
  &\leq\int_{\Omega}\tdens_0\dx+(\Lyap(\tdens_0))^{\frac{1}{2}}
\end{align*}

\subsection{Proof of Proposition~\ref{u:C1}: $\mathcal{C}^1$-regularity of $u(t)$}
\label{sec:proof-regularity}

Let $t\in[0,\tau(\tdens_0)[$ and choose $h>0$ such that
$t+h<\tau(\tdens_0)$. 
The solution $\tdens(t)$ of~\eqref{CP} belongs to a ball
$B(\tdens_0,R)$ centered in $\tdens_0$ and with appropriate radius $R$. 
More detailed information on the size of this $R$ can be drawn from
the proof of Theorem~\ref{exist:ODE}. 
Using~\eqref{always:pos}, it is possible to find two constants
$a(\tdens_0),b(\tdens_0)$ such that
$\tdens(t)\in\mathcal{D}\left(a(\tdens_0),b(\tdens_0)\right)$. 
This allows us to write:
\begin{equation*}
  \Kmu(\tdens(t))=
     \frac{1}{\lambda(\tdens(t))}
     \left(
       \frac{\|\tdens(t)\|_{\Cac}}{\lambda(\tdens(t))}
     \right)^{\frac{n+\delta}{2\delta}}
  \leq 
      \frac{1}{a(\tdens_0)}
      \left(
        \frac{b(\tdens_0)}{a(\tdens_0)}
      \right)^{\frac{{n+\delta}}{2\delta}}
   \quad  \forall t\in[0,\tau(\tdens)[
\end{equation*}
which shows that both the Potential and Flux operators are bounded and
Lipschitz-continuous in
$\mathcal{D}\left(a(\tdens_0),b(\tdens_0)\right)$. 

Heuristically, the proof is based on the observation that, assuming
that both $\tdens'$ and $u'$ exist, we can take the derivative in time
of~\eqref{sys:div:weak} and~\eqref{bou:cond:u:weak} and use
the fact that the source function is independent on time, thus
obtaining~\eqref{u'}.
We first note that $u(t)=\Pot(\tdens(t))$ is Lipschitz-continuous
in time, since $\Pot$ is locally Lipschitz-continuous and
$\tdens\in\mathcal{C}^1\left([0,\tau(\tdens_0)[\,;\Cac\right)$. 
Next at each time $t\in[0,\tau(\tdens_0)[$, we define $w(t)$, that 
heuristically should be $u'(t)$, as the unique solution of:
\begin{equation}
  \label{w}
  \left\{
    \begin{array}{cc}
      \int_{\Omega}\tdens(t)\nabla w(t)\cdot\nabla\varphi\dx =
        -\int_{\Omega}\tdens'(t)\nabla u(t)\cdot\nabla\varphi\dx
      \quad \forall \varphi \in H^1(\Omega) \\
   \int_{\Omega} w(t)\dx=0
   \end{array}
   \right.
\end{equation}
It is easy to verify that $w(t)\in\Hac$.
Hence $\forall \varphi \in H^1(\Omega)$ we can write:
\begin{equation*}
  \int_{\Omega}\tdens(t)\nabla u(t)\cdot\nabla\varphi\dx=
      \int_{\Omega} f\varphi\dx=
       \int_{\Omega}\tdens(t+h)\nabla u(t+h)\cdot\nabla\varphi\dx 
\end{equation*}
Changing sign and adding in both sides 
$\int_{\Omega}\tdens(t)\nabla u(t+h)\cdot\nabla\varphi\dx$, yields:
\begin{multline}
  \label{diff:d:2} 
    \int_{\Omega}\tdens(t)\nabla [u(t+h)-u(t)] \cdot \nabla \varphi\dx= \\
    -\int_{\Omega}[\tdens(t+h)-\tdens(t)]\nabla
       u(t+h)\cdot\nabla\varphi\dx
       \qquad
\end{multline}
Now we multiply~\eqref{w} by $-h$ and sum~\eqref{w}
and~\eqref{diff:d:2} to obtain:
\begin{equation*}
  \begin{aligned}
  \int_{\Omega}\tdens(t)
  &    \nabla[u(t+h)-u(t)-hw(t)]\cdot\nabla\varphi\dx\\
  &=-\int_{\Omega}[\tdens(t+h)-\tdens(t)] 
      \nabla u(t+h)\cdot\nabla\varphi\dx
   +h\int_{\Omega}\tdens'(t)\nabla u(t)\cdot\nabla\varphi\dx \\
  &= -\int_{\Omega}
     \left\{\left[\tdens(t+h)-\tdens(t)-h\tdens'(t)\right]\nabla
       u(t+h)\right.\\ 
  &\qquad\qquad\qquad\qquad
   +\left.h\tdens'(t)\left(\nabla u(t+h)-\nabla u(t)\right)\right\}
     \cdot\nabla\varphi\dx \\
  &=-\int_{\Omega}[G_1(t,h)+h\; G_2(t,h)]\cdot \nabla\varphi\dx
  \end{aligned}
\end{equation*}
with 
\begin{align*}
  G_1(t,h)&=\left[\tdens(t+h)-\tdens(t)-h\tdens'(t)\right]\nabla u(t+h)\\
  G_2(t,h)&=\tdens'(t)\left[\nabla u(t+h)-\nabla u(t)\right]
\end{align*}
Since $\tdens\in \mathcal{C}^1(0,\tau;\Cac)$, we can estimate the
above functions $G_1$ and $G_2$ as:
\begin{align*}
  \|G_1(t,h)\|_{\Cac}
    &\leq \|\tdens(t+h)-\tdens(t)-h\tdens'(t)\|_{\Cac}\|\nabla u(t+h)\|_{\Cac}\\
    &\leq
      \Kdelta(n,\Omega,\delta)\Kmu(\tdens(t))\|f\|_{L^{\infty}(\Omega)}
      \cdot o(h)\\
    &\leq
      \Kdelta(n,\Omega,\delta)K(\tdens_0)\|f\|_{L^{\infty}(\Omega)}
      \cdot o(h)
\end{align*}
and, since the Potential operator is Lipschitz-continuous, we have
also: 
\begin{align*}
  \|G_2(t,h)\|_{\Cac}
  &=\|\tdens'(t)\left[\nabla u(t+h)-\nabla u(t)\right]\|_{\Cac}\\
  &\leq \| \tdens'(t)\|_{\Cac}\|\nabla u(t+h)-\nabla u(t)\|_{\Cac}\\
  &\leq L(\tdens_0)h
\end{align*}
where $L(\tdens_0)$ is a function of $f$, $\Kdelta$.
Thus we can write:
\begin{equation*}
  \|G_1+hG_2\|_{\Cac}\leq \|G_1\|_{\Cac}+\|G_2\|_{\Cac} = o(h)
\end{equation*}
and, for Theorem~\ref{holder:f:G} using $F=-(G_1+hG_2)$, we obtain:
\begin{equation*}
\lim_{h\rightarrow 0} \frac{\|\nabla [u(t+h)-u(t)-hw(t)]\|_{\Cac}}{h}=0
\end{equation*}
that shows that $u'\in \mathcal{C}^{1}(0,\tau;\Cac)$ with $u'=w$.

\subsection{Proof of Lemma~\ref{holder:f:G}: elliptic regularity}
\label{proof:holder:f:G}

Lemma~\ref{holder:f:G} is analogous to Theorem~5.19
of~\citet{giaquinta} simplified to a scalar elliptic equation but
extended to explicitly determine the dependence of the inequality
constants upon $\tdens$. 
We will denote with $C$ or $c$ generic constants that may depend upon
$n$, $\Omega$, and the H\"older continuity exponent $\delta$ but are
always independent of $\tdens$. 
We will use the following result adapted from Proposition 5.8
in~\citet{giaquinta}: 
\begin{Lemma}[Elliptic Decay]
  \label{decay}
  Let $v\in H^1(\Omega)$ be any solution of 
  \begin{equation}
    \label{laplace}
    \int_{\Omega}\nabla v\nabla \varphi\dx=0 
       \quad \forall \varphi\in H^1_0(\Omega)
  \end{equation}
  then there exists a constant $c(n)$ such that:
  \begin{gather}
    \label{estimate:n}
    \int_{\Ball[\rho]{x_0}}|\nabla v|^2\dx
    \leq c(n)\left(\frac{\rho}{R}\right)^n \int_{\Ball[R]{x_0}}|\nabla v|^2 \dx\\
    \label{estimate:n+2}
    \int_{\Ball[\rho]{x_0}}|\nabla v-\Avgint[x_0]{\rho}{\nabla v}|^2\dx
    \leq c(n)\left(\frac{\rho}{R}\right)^{n+2} 
        \int_{\Ball[R]{x_0}}|\nabla v-\Avgint[x_0]{R}{\nabla v}|^2 \dx
  \end{gather}
  for arbitrary balls $\Ball[\rho]{x_0}\Subset \Ball[R]{x_0}\Subset\Omega$.
\end{Lemma}

\noindent
Our case follows from the observation that the derivatives of $v$
satisfy the weak form of Laplace equation 
(see also~\citet{ambrosio_ln}, page 61).
Note that the constant $c(n)$ depends only on the problem dimension
$n$ as we are considering Laplace equation.

We also use the following result from Lemma 5.13 in~\citet{giaquinta}
and Lemma 9.2 in~\citet{Ambrosio:2004}:
\begin{Lemma}[Iteration lemma]
  \label{iteration:lemma}
  Let $\phi:\REAL^+\mapsto \REAL^+$ be a non-negative and non
  increasing function satisfying 
  \begin{equation}
    \label{hp:lemma}
    \phi(\rho)\leq A\left[\left(\frac{\rho}{R}\right)^{\alpha}
      +\epsilon\right]\phi(R)+B\; R^\beta
  \end{equation}
  for some $A,\alpha,\beta>0$, with $\alpha>\beta$ and for all
  $0<\rho\leq R\leq R_0$, where $R_0>0$ is given. 

  Then there exist constants $\epsilon_0=\epsilon_0(A,\alpha,\beta)$
  and $C=C(A,\alpha,\beta)$ such that
  \begin{equation}
    \label{concl:lemma}
    \text{if }\quad
    \epsilon\leq \epsilon_0=\left(\frac{1}{2A}\right)^{\frac{2\alpha}{\alpha-\beta}} 
    \quad \text{ then } \quad
    \phi(\rho)\leq C\left[\frac{\phi(R)}{R^{\beta}}+B \right] \rho^\beta.
  \end{equation}
\end{Lemma}

The proof of Lemma~\ref{holder:f:G} uses a classical bootstrap
technique introduced by~\citet{Morrey:1954,Campanato:1965} and used
more recently by~\citet{ColomboMingione:2015} to show the regularity
of local minimizers of double phase variational integrals. The
technique can be described by the following steps.  First we consider
a compact set $K\Subset\Omega$ and prove that $u\in \mor{\nu}{K}$ for
a suitable regularity exponent $\nu$ with $0<\nu<n$. Then,
$u\in \cam{n+2\delta}{K}$ where $\mor{\nu}{K} $ and
$\cam{n+2\delta}{K}$ are the Morrey and Campanato spaces,
respectively.  The results are extended to the entire domain by
assuming enough regularity of $\partial\Omega$.  This latter steps is
not reported in the following proof for brevity.  Finally, the
equivalence between the Campanato spaces $\cam{n+2\delta}{\Omega}$ and
$\Cac$ is used to prove estimate~\eqref{grad:u:Cac} and to derive the
expression of the constant $\Kmu$ given in~\eqref{estimate:C(d)}.  We
recall that the norm of a function $u:\Omega\rightarrow\REAL^m$ (in
our case we have either $m=1$ or $m=n$) belonging to a Morrey space is
given by:
\begin{equation*}
  \|u\|_{\mor{\gamma}{\Omega}} = 
     \Big(\sup_{\underset{\rho>0}{x_0\in\Omega}}
     \rho^{-\gamma}\int_{\Omega(x_0,\rho)} |u|^2 \dx\Big)^{\frac{1}{2}}
\end{equation*}
where $\Omega(x_0,\rho)=\Omega\cap \Ball[\rho]{x_0}$ and $0\leq\gamma<n$.
For $0\leq\gamma<n+2$, the norm of $u$ belonging to a Campanato space
is given by:
\begin{equation*}
  \|u\|_{\cam{\gamma}{\Omega}} = 
     \|u\|_{L^2(\Omega)}
     + \Big(\sup_{\underset{\rho>0}{x_0\in\Omega}}
     \rho^{-\gamma} \int_{\Omega(x_0,\rho)}
     |u-\Avgint[x_0]{\rho}{u}|^2 
      \dx\Big)^{\frac{1}{2}}
\end{equation*}
where $\Avgint[x_0]{\rho}{u}=\int_{\Omega(x_0,\rho)} u\dx/|\Omega(x_0,\rho)|$ 
is the average integral.

\begin{proof}
The first step of the bootstrap proceeds as follows.
Consider $x_0\in K$ and the ball $\Ball[R]{}:=\Ball[R]{x_0}\Subset\Omega$. 
In this ball we use Korn's technique (freezing the coefficients) to
decompose the solution as $u=v+w$ where $v\in H^1(\Ball[R]{})$ satisfies
the equations: 
\begin{equation}
  \label{eq:v}
   \int_{\Ball[R]{}} \tdens(x_0) \nabla v \nabla \varphi\dx=0 
         \quad \forall \varphi \in H^1_0(\Ball[R]{})
\end{equation}
with $v=u$ in $\partial \Ball[R]{}$ and the second equation is to be
interpreted in the sense that 
$v-u\in H^1_0(\Ball[R]{})$. 
The second function $w\in H^1_0(\Ball[R]{})$ satisfies the equation:
\begin{align}
  \label{eq:w}
  &
     \int_{\Ball[R]{}} \tdens(x_0) \nabla w \nabla \varphi\dx= 
    \int_{\Ball[R]{}} \Big[F_0\varphi+F\cdot\nabla\varphi-
    \nonumber \\
  &\qquad\qquad\qquad\qquad
    \left(\tdens(x)-\tdens(x_0)\right)
    \nabla u\cdot\nabla\varphi\Big]\dx \qquad
    \forall \varphi \in H^1_0(\Ball[R]{}) 
\end{align}
with $w=0$ in $\partial \Ball[R]{}$.
Since $\tdens(x_0)$ in~\eqref{eq:v} is a strictly positive and
bounded scalar number it can be eliminated from the equation, hence
$w$ simply solves the weak form of Laplace equation:
\begin{equation}
   \int_{\Ball[R]{}} \nabla v \nabla \varphi\dx=0 
      \qquad \forall \varphi \in H^1_0(\Omega)
  \label{eq:v:laplace}
\end{equation}
with $v=u$ in $\partial\Omega$.
Thus we can use Lemma~\ref{decay} to obtain:
\begin{equation}
  \label{estimate:laplace}
  \int_{\Ball[\rho]{}}|\nabla v|^2\dx\leq 
    c(n)\left(\frac{\rho}{R}\right)^n \int_{\Ball[R]{}}|\nabla v|^2\dx
\end{equation}
Recall that at this point our goal is to estimate the Morrey norm 
$\|\nabla u\|_{\mor{\nu}{K}}$ with $\nu<n$.
We use the above decomposition of $u$ to estimate
$\phi(\rho):=\int_{\Ball[\rho]{}}|\nabla u|^2\dx$, $0<\rho\leq R$. 
Thus we can write:
\begin{align*}
  \int_{\Ball[\rho]{}}|\nabla u|^2\dx
  & =\int_{\Ball[\rho]{}}|\nabla v+\nabla w|^2\dx
    \leq
    2\int_{\Ball[\rho]{}}|\nabla v|^2\dx
    + 2\int_{\Ball[\rho]{}}|\nabla w|^2\dx \\
  & \leq 
    c(n)\left(\frac{\rho}{R}\right)^n \int_{\Ball[R]{}}|\nabla v|^2\dx +
    2\int_{\Ball[\rho]{}}|\nabla w|^2\dx \\
  & = c(n)\left(\frac{\rho}{R}\right)^n 
    \int_{\Ball[R]{}}|\nabla u-\nabla w|^2\dx
    + 2\int_{\Ball[\rho]{}}|\nabla w|^2\dx \\
  & \leq c(n)\left(\frac{\rho}{R}\right)^n
        \int_{\Ball[R]{}}|\nabla u|^2\dx
     + c(n)\left(\frac{\rho}{R}\right)^n\int_{\Ball[\rho]{}}|\nabla w|^2\dx \\
  & \qquad   + 2\int_{\Ball[\rho]{}}|\nabla w|^2\dx \\
  & \leq c(n)\left(\frac{\rho}{R}\right)^n
    \int_{\Ball[R]{}}|\nabla u|^2\dx
    +c(n)\int_{\Ball[R]{}}|\nabla w|^2\dx
\end{align*}
Note that, somewhat improperly, we always use the symbol $c(n)$
to indicate a constant depending on $n$ only and that may assume
different meaning even within the same equation. 
To estimate $\displaystyle{\int_{\Ball[R]{}}|\nabla w|^2\dx}$ we use
$\varphi=w$ in~\eqref{eq:w} to get:
\begin{equation}
  \label{pre:ine}
  \begin{split}
    \lambda(\tdens)\int_{\Ball[R]{}}|\nabla w|^2\dx
    & \leq \int_{\Ball[R]{}}\tdens(x_0)|\nabla w|^2\dx  \\
    & =\int_{\Ball[R]{}}\left[
      F_0w+F\cdot\nabla w-\left(\tdens(x)-\tdens(x_0)\right)
      \nabla u\nabla w\right]\dx
  \end{split}
\end{equation}
Using  H\"older continuity of $\tdens$, and Poincar\'e and
Cauchy-Scharwz inequalities, we can bound the
right-hand-side of the previous equation to obtain: 
\begin{align*}
  \int_{B(x_0,R)} F_0 w\dx 
  & \leq 
  \|F_0\|_{L^2(B(x_0,R))}\;c(n)\,\|\nabla w\|_{L^2(B(x_0,R))}\\
  \int_{B(x_0,R)} F\cdot\nabla w\dx
  & \leq
  \|F\|_{L^2(B(x_0,R))}\;\|\nabla w\|_{L^2(B(x_0,R))}\\
  \int_{B(x_0,R)} \left(\tdens(x)-\tdens(x_0)\right)\nabla u\nabla w \dx
  & \leq 
  R^{\delta} \|\tdens\|_{\Cac}\|\nabla u\|_{L^2(B(x_0,R))}\,\|\nabla w\|_{L^2(B(x_0,R))}
\end{align*}
In the end, using Minkowski inequality to remove the double products,
we can write: 
\begin{equation}
  \label{lax:migram:w}
  \begin{split}
    \int_{\Ball[R]{}}|\nabla w|^2\dx
    \leq 
    2\frac{1}{\left(\lambda(\tdens)\right)^2}
    \left[
      (c(n))^2\|F_0\|^2_{L^2(\Ball[R]{})} \right.\\
      \left. +\|F\|^2_{L^2(\Ball[R]{})}+
      R^{2\delta} \|\tdens\|^2_{\Cac} \|\nabla u\|^2_{L^2(\Ball[R]{})}
    \right]
  \end{split}
\end{equation}
Since $F_0\in L^{\infty}(\Omega)$, implying that
$F_0\mor{\nu}{\Omega}$ and
$\|F_0\|^2_{\mor{\nu}{\Omega}}\leq c(n)\|F_0\|^2_{L^{\infty}(\Omega)}$,
we obtain: 
\begin{equation}
  \label{f:L:infty}
  \|F_0\|^2_{L^2(\Ball[R]{})}
  \leq c(n)\|F_0\|^2_{L^{\infty}(\Ball[R]{})}R^\nu
  \leq c(n)\|F_0\|^2_{L^{\infty}(\Omega)}R^\nu 
\end{equation}
Since $F\in\Cac$, each component of 
$F$ belongs to $\cam{\gamma}{\Omega}$ for all $0\leq\gamma\leq n+2\delta$. 
Noting that we require $0\leq\nu<n$,
and in this case $\mor{\nu}{\Omega}\equiv\cam{\nu}{\Omega}$,
we obtain:
\begin{equation}
  \label{G:Cac:mor}
  \|F\|^2_{L^2(\Ball[R]{})}
  \leq\|F\|^2_{\cam{\gamma}{\Ball[R]{}}} R^\gamma
  \leq c(n)\|F\|^2_{\Cac}R^\gamma
\end{equation}
Taking $\nu<n$ in~\eqref{f:L:infty} 
and $\gamma=\nu$  in~\eqref{G:Cac:mor} 
we get:
\begin{align}
  \label{pre:lemma}
  \int_{\Ball[\rho]{}}|\nabla u|^2\dx
  & \leq c(n)
    \left[
       \left(\frac{\rho}{R}\right)^n+ 
       R^{2\delta}\left(\frac{\|\tdens\|_{\Cac}}{\lambda(\tdens)}\right)^2
    \right]
    \int_{\Ball[R]{}}|\nabla u|^2\dx \nonumber\\ 
  & \qquad +c(n)
    \left( 
       \frac{
           \|F_0\|^2_{L^\infty(\Omega)}+\|F\|^2_{\Cac}
        }{
           \left(\lambda(\tdens)\right)^2
    }\right)R^{\nu}
\end{align}
Now we rewrite inequality~\eqref{pre:lemma} 
in the form of the hypotheses of Lemma~\ref{iteration:lemma}, i.e.:
\begin{gather*}
  \phi(\rho):=\int_{\Ball[\rho]{}}|\nabla u|^2\dx,\qquad \alpha=n,\  \beta=\nu, \\ 
  \epsilon=R^{2\delta}  \left(\frac{\|\tdens\|_{\Cac} }{\lambda(\tdens)}\right)^2,
  \;A=c(n),\;
  B=c(n)
  \left( 
    \frac{\|F_0\|^2_{L^\infty(\Omega)}+\|F\|^2_{\Cac}}{(\lambda(\tdens))^2}
  \right)
\end{gather*}
for $\rho\leq R$. Considering $R$ such that: 
\begin{equation*}
  R^{2\delta} \left(\frac{\|\tdens\|_{\Cac}}{\lambda(\tdens)}\right)^2 
  \leq \left(\frac{1}{2A}\right)^{\frac{2n}{n-\nu}}=A_0
\end{equation*}
we have that:
\begin{equation}
  R\leq R_0=A^{\frac{1}{\delta}}_0
     \left(
         \frac{\lambda(\tdens)}{\|\tdens\|_{\Cac}}
     \right)^{\frac{1}{\delta}}
     \label{R_0}
\end{equation}
We can now apply Lemma~\ref{iteration:lemma} to obtain,
for $0<\rho\leq R\leq R_0$, the following estimate:
\begin{equation}
  \label{after:lemma}
  \int_{\Ball[\rho]{}}|\nabla u|^2\dx
  \leq C(A,n,\nu)\rho^\nu 
  \left(\frac{\int_{\Ball[R]{}}|\nabla u|^2\dx}{R^{\nu}}+B\right)
\end{equation}
Incorporating all the constant into one single constant $C(n,\nu)$ we
obtain: 
\begin{equation}
  \label{morrey:inizio}
  \int_{\Ball[\rho]{}}|\nabla u|^2\dx
  \leq C(n,\nu)\rho^\nu 
  \left(
    \frac{\int_{\Ball[R]{}}|\nabla u|^2\dx}{R^{\nu}}
    + \frac{\|F_0\|^2_{L^\infty(\Omega)}+\|F\|^2_{\Cac}}{(\lambda(\tdens))^2}
  \right)
\end{equation}
The previous estimate is valid for every $\Ball[R]{}\Subset\Omega$. 
Varying $x_0\in K$ and using the continuity inequality in the
Lax-Milgram Lemma we obtain the desired estimate of this first step of
the bootstrap procedure, i.e.:
\begin{align}
 \label{morrey:intern}
 \|\nabla u\|^2_{\mor{\nu}{K}}
  &\leq C(n,\nu)\left(
      \frac{\int_{\Ball[R_0]{}}|\nabla u|^2\dx}{R_0^{\nu}}+
      \frac{\|F_0\|^2_{L^\infty(\Omega)}+\|F\|^2_{\Cac}}{(\lambda(\tdens))^2}
    \right) \nonumber\\
 & \leq C(n,\nu) \left(
       \frac{
          C(\Omega)\|F_0\|^2_{L^\infty(\Omega)}+\|F\|^2_{\Cac}
        }{
          \left(\lambda(\tdens)\right)^2
        }
        \frac{1}{R_0^{\mu}}+
        \frac{\|F_0\|^2_{L^\infty(\Omega)}+\|F\|^2_{\Cac}}{(\lambda(\tdens))^2}
   \right)\nonumber\\
 &\leq C(n,\nu)C(\Omega)
        \frac{\|F_0\|^2_{L^\infty(\Omega)}+\|F\|^2_{\Cac}}{(\lambda(\tdens))^2}
        \left(
          \left(
            \frac{\|\tdens\|_{\Cac}}{\lambda(\tdens)}
           \right)^{\frac{\nu}{\delta}} + 1
        \right) \nonumber\\
  &\leq C(n,\nu)C(\Omega)
    \frac{\|F_0\|^2_{L^\infty(\Omega)}+\|F\|^2_{\Cac}}{(\lambda(\tdens))^2}
    \left(\frac{\|\tdens\|_{\Cac}}{\lambda(\tdens)}\right)^{\frac{\nu}{\delta}}
\end{align}
where $C(n,\nu)$ is bounded for all $\nu<n$.

The second step of the bootstrap procedure starts by noting that
Equation~\eqref{eq:w} can be rewritten using $R=R_0$ as defined above:
\begin{equation}
  \begin{split}
    & 
    \int_{\Ball[R]{}}\tdens(x_0)\nabla w\nabla\varphi\dx= 
    \int_{\Ball[R]{}}
    \Big[F_0\varphi+\left(F-\Avgint{R}{F}\right)\cdot\nabla\varphi-
    \\
    & 
    \qquad\qquad 
      \left(\tdens(x)- \tdens(x_0)\right)\nabla u\cdot\nabla\varphi 
    \Big]\dx 
    \qquad
    \forall \varphi \in H^1_0(\Ball[R]{})\\
    & w=0 \text{ in } \ \partial \Ball[R]{}
  \end{split}
  \label{eq:w:bis}
\end{equation}
We continue by using again the decomposition $u=v+w$ and
Lemma~\ref{decay} to obtain:
\begin{align*}
  \int&_{\Ball[\rho]{}}|\nabla u-\Avgint{\rho}{\nabla u}|^2\dx
   = 
    \int_{\Ball[\rho]{}}|\nabla v+\Avgint{\rho}{\nabla v}
    +\nabla w+\Avgint{\rho}{\nabla w}|^2\dx \\
  & \leq  
    c(n)\left(\frac{\rho}{R}\right)^{n+2} 
    \int_{\Ball[R]{}}|\nabla v-\Avgint{R}{\nabla v}|^2\dx 
   +2\int_{\Ball[\rho]{}}|\nabla w-\Avgint{\rho}{\nabla w}|^2 \dx \\
  &\leq   c(n)\left(\frac{\rho}{R}\right)^{n+2}
    \int_{\Ball[R]{}}|\nabla u-\Avgint{R}{\nabla u}|^2\dx
    +c(n)\int_{\Ball[R]{}}|\nabla w|^2\dx
\end{align*}
where the last inequality arises from the minimality of the mean.
We follow the same developments as before, but now we explicitly
include the factor $R$ in the constant of Poincar\'e inequality to
obtain: 
\begin{gather}
  \notag
  \int_{\Ball[\rho]{}}|\nabla u-\Avgint{\rho}{\nabla u}|^2\dx \leq
  c(n)\left(\frac{\rho}{R}\right)^{n+2}
  \int_{\Ball[R]{}}|\nabla u-\Avgint{R}{\nabla u}|^2\dx \\
  + 2c(n)
  \frac{
    R^2\|F_0\|^2_{L^2(\Ball[R]{})}+\|F-\Avgint{R}{F}\|_{L^2(\Ball[R]{})}+
    R^{2\delta}\|\tdens\|^2_{\Cac}\|\nabla u\|^2_{L^2(\Ball[R]{})}
  }{
    \left(\lambda(\tdens)\right)^2
  }
 \label{pre:mor}
\end{gather}
Since $\nabla u\in\mor{\nu}{K}$ for $0<\nu<n$ we can take $\nu=n-\delta$ in  
\eqref{morrey:intern} to get:
\begin{equation*}
  \|\nabla u\|^2_{L^2(\Ball[R]{})}=
  \frac{\int_{B_{R}}|\nabla u|^2\dx}{R^{n-\delta}}
  R^{n-\delta}\leq\|\nabla u\|_{L^{2,n-\delta}}R^{n-\delta}
\end{equation*}
Using $\nu=n-2+\delta$ in~\eqref{f:L:infty} we obtain
  $\|F_0\|^2_{L^2(\Ball[R]{})}\leq c(n)\|F_0\|_{L^{\infty}(\Omega)}R^{n-2+\delta}$,
%
while using $\gamma=n+\delta$ in~\eqref{G:Cac:mor} we have
  $\|F-\Avgint{R}{F}\|_{L^2(\Ball[R]{})}\leq \|F\|_{\Cac}R^{n+\delta}$.
%
Substitution of these inequalities in \eqref{pre:mor} yields:
\begin{align*}
  &\int_{\Ball[\rho]{}}|\nabla u-\Avgint{\rho}{\nabla u}|^2\dx
  \leq c(n)
    \left(\frac{\rho}{R}\right)^{n+2}
    \int_{\Ball[R]{}}|\nabla u-\Avgint{R}{\nabla u}|^2\dx\\
  &\quad
    + c(n)
    \frac{\|F_0\|^2_{L^\infty(\Omega)}+\|F\|_{\Cac}}{(\lambda(\tdens))^2}R^{n+\delta}\\
  &\quad
    +R^{2\delta} \frac{\|\tdens\|^2_{\Cac}}{(\lambda(\tdens))^2} 
    C(n,n-\delta) C(\Omega)
    \frac{\|F_0\|^2_{L^\infty(\Omega)}+\|F\|^2_{\Cac}}{(\lambda(\tdens))^2}
    \left(
       \left(
          \frac{\|\tdens\|_{\Cac}}{\lambda(\tdens)}
       \right)^{\frac{n-\delta}{\delta}}
    \right)R^{n-\delta}\\
 & \leq c(n)\left(\frac{\rho}{R}\right)^{n+2}
   \int_{\Ball[R]{}}|\nabla u-\Avgint{R}{\nabla u}|^2\dx\\
 &\quad
   +  C(n,\Omega,\delta)
   \frac{\|F_0\|^2_{L^\infty(\Omega)}+\|F\|^2_{\Cac}}{(\lambda(\tdens))^2}
   \left(
      1+\frac{\|\tdens\|^2_{\Cac}}{(\lambda(\tdens))^2}
      \left(
        \frac{\|\tdens\|_{\Cac}}{\lambda(\tdens)}
      \right)^{\frac{n-\delta}{\delta}}
   \right)R^{n+\delta}\\
 & \leq c(n)
   \left(\frac{\rho}{R}\right)^{n+2}
      \int_{\Ball[R]{}}|\nabla u-\Avgint{R}{\nabla u}|^2\dx\\
 &\quad
   + C(n,\Omega,\delta)
   \frac{\|F_0\|^2_{L^\infty(\Omega)}+\|F\|^2_{\Cac}}{(\lambda(\tdens))^2}
   \left(
      \frac{\|\tdens\|_{\Cac}}{\lambda(\tdens)}
   \right)^{\frac{n+\delta}{\delta}} R^{n+\delta}
\end{align*}
Applying Lemma~\ref{iteration:lemma} with 
$\phi(\rho):=\int_{\Ball[\rho]{}}|\nabla u-\Avgint{\rho}{\nabla u}|^2\dx$
yields for $0<\rho\leq R\leq R_0$:
\begin{gather*}
  \int_{\Ball[\rho]{}}|\nabla u-\Avgint{\rho}{\nabla u}|^2\dx \leq 
   \rho^{n+\delta}C(n,\Omega,\delta) \cdot\\
  \cdot \left[\frac{\int_{\Ball[R]{}}|\nabla u-\Avgint{R}{\nabla u}|^2\dx}{R^{n+\delta}}+
    \left(\frac{\|F_0\|^2_{L^\infty(\Omega)}+\|F\|^2_{\Cac}}{\lambda(\tdens)^2}\right)
    \left(\frac{\|\tdens\|_{\Cac}}{\lambda(\tdens)}\right)^{\frac{n}{\delta}+1}\right]
\end{gather*}
from which, using again the minimality of the mean and the
estimate of $R_0$ given in~\eqref{R_0}, we can evaluate:
\begin{align*}
  &\frac{
       \int_{\Ball[\rho]{}}|\nabla u-\Avgint{\rho}{\nabla u}|^2\dx
    }{
       \rho^{n+\delta}
    } \leq C(n,\delta,\Omega) \cdot \\
  &\quad 
    \cdot    \left[
      \frac{\int_{\Ball[R_0]{}}|\nabla u-\Avgint{R_0}{\nabla u}|^2\dx}{R_0^{n+\delta}}+
      \left(
        \frac{\|F_0\|^2_{L^\infty(\Omega)}+\|F\|^2_{\Cac}}{\lambda(\tdens)^2}
      \right)
      \left(\frac{\|\tdens\|_{\Cac}}{\lambda(\tdens)}\right)^{\frac{n}{\delta}+1}
    \right]\\
  &\quad \leq  C(n,\delta,\Omega)
    \left(
      \int_{\Ball[R_0]{}}\!\!\!|\nabla u-\Avgint{R_0}{\nabla u}|^2\dx+
        \frac{\|F_0\|^2_{L^\infty(\Omega)}+\|F\|^2_{\Cac}}{\lambda(\tdens)^2}
    \right)
    \left(
      \frac{\|\tdens\|_{\Cac}}{\lambda(\tdens)} 
    \right)^{\frac{n}{\delta}+1}\\
  &\quad \leq  C(n,\delta,\Omega)
    \left( 
       \int_{\Omega}|\nabla u|^2\dx+
         \frac{\|F_0\|^2_{L^\infty(\Omega)}+\|F\|^2_{\Cac}}{\lambda(\tdens)^2} 
    \right)
     \left(
         \frac{\|\tdens\|_{\Cac}}{\lambda(\tdens)}
     \right)^{\frac{n}{\delta}+1}\\
  &\quad \leq  C(n,\delta,\Omega)
    \left(
       \frac{\|F_0\|^2_{L^\infty(\Omega)}+\|F\|^2_{\Cac}}{\lambda(\tdens)^2}
    \right)
    \left(\frac{\|\tdens\|_{\Cac}}{\lambda(\tdens)}\right)^{\frac{n}{\delta}+1}
\end{align*}
Hence $\nabla u\in\cam{n+\delta}{K}$ and we can write:
\begin{equation}
  \|\nabla u\|^2_{\cam{n+\delta}{K}}\leq 
  C(n,\delta,\Omega)
  \left(\frac{\|F_0\|^2_{L^\infty(\Omega)}+\|F\|^2_{\Cac}}{\lambda(\tdens)^2}\right)
  \left(\frac{\|\tdens\|_{\Cac}}{\lambda(\tdens)}\right)^{\frac{n}{\delta}+1}
  \label{second:step}
\end{equation}
The bootstrap procedure is restarted from~\eqref{eq:w:bis}
using $\nu=n-2+2\delta$ in~\eqref{f:L:infty}  
and $\gamma=n+2\delta$ in~\eqref{G:Cac:mor},
and estimate~\eqref{second:step} in~\eqref{pre:mor} so that a term
$R^{n+2\delta}$ can be factored. 
Thus we can write:
\begin{align*}
  \int_{\Ball[\rho]{}}|\nabla u&-\Avgint{\rho}{\nabla u}|^2\dx
   \leq c(n)
    \left(\frac{\rho}{R}\right)^{n+2}
       \int_{\Ball[R]{}}|\nabla u-\Avgint{R}{\nabla u}|^2\dx \\
  & + c(n)
    \frac{\|F_0\|^2_{L^\infty(\Omega)}+\|F\|_{\Cac}}{(\lambda(\tdens))^2}R^{n+2\delta}\\
  & + C(n,\delta,\Omega)
    \frac{\|F_0\|^2_{L^\infty(\Omega)}+\|F\|^2_{\Cac}}{\lambda(\tdens)^2}
    \left(\frac{\|\tdens\|_{\Cac}}{\lambda(\tdens)}\right)^{\frac{n}{\delta}+1}
    R^{n+2\delta}
\end{align*}
and finally, applying once again lemma \ref{iteration:lemma}, we have
the final result:
\begin{equation}
  \label{estimate:cam}
  \|\nabla u\|^2_{\cam{n+2\delta}{K}} \leq
  C(n,\delta,\Omega)
  \frac{\|f\|^2_{L^\infty(\Omega)}+\|G\|^2_{\Cac}}{\lambda(\tdens)^2}
  \left(\frac{\|\tdens\|_{\Cac}}{\lambda(\tdens)}\right)^{\frac{n}{\delta}+1}
\end{equation}
Extension of the previous estimate to the entire domain $\Omega$ can
be obtained following the same bootstrap procedure starting from the
analogue of the elliptic decay lemma~\ref{decay} on hemispheres
(similarly to what is proposed in~\citet{giaquinta}, Theorem 5.21).
Such process introduces a dependence on the regularity of the boundary
$\partial\Omega$ in the constant $C(n,\delta,\Omega)$
in~\eqref{estimate:cam}, but we do not explicitly write such
dependence. 
By the equivalence between $\cam{n+2\delta}{\Omega}$ and $\Cac$ 
we get:
\begin{equation*}
  \|\nabla u\|_{\Cac}\leq C(n,\delta,\Omega)
     \frac{\|F_0\|_{L^\infty(\Omega)}+\|F\|_{\Cac}}{\lambda(\tdens)}
     \left(\frac{\|\tdens\|_{\Cac}}{\lambda(\tdens)}\right)^{\frac{n+\delta}{2\delta}}
\end{equation*}
which proves~\eqref{grad:u:Cac} and~\eqref{estimate:C(d)}.
From this, using Theorem 1.40 of~\citet{troi}, we directly obtain that
$u\in\mathcal{C}^{1,\delta}(\bar{\Omega})$. 
\end{proof}

\section{Conclusions}

We have presented a continuous extension or the original model
governing the dynamics of the Physarum Polycephalum slime mold. 
The proposed model couples an elliptic diffusion equation enforcing
PP density balance with an ordinary differential equations governing the
dynamics of the flow of information along the PP body.
By transforming the problem into an infinite dimensional dynamical
system, we are able to prove local in time existence and uniqueness of
the problem solution under the hypothesis $L^\infty(\Omega)$ forcing
function and H\"older continuous diffusion coefficient $\tdens$.
Extension to larger (possibly infinite) times is not currently at
reach, but the existence of a Lyapunov-candidate function and
numerical evidence more seems to strongly justify our belief that a
steady state solution to our model problem exists and that at this
steady state conditions our model problem is equivalent to the
PDE-based Monge-Kantorovich Optimal Transport problem. 
The Lyapunov-candidate function, which is characterized by a
non-positive Lie derivative, is derived by recalling the analogy between
the Monge-Kantorovich equations and shape optimization problems. 

We describe a simple but rather efficient numerical formulation for
the quantitative solution of our model problem and the simulation of
the dynamics of Physarum Polycephalum.
The numerical approach is based on 
a combination of piecewise constant a FEM discretization space for
$\tdens$ and a piecewise linear FEM space, defined on a once-refined
triangulation, for the potential $u$. 
A simple forward Euler time discretization completes the discrete
formulation. 
Preliminary numerical results show that the dynamics of the slime-mold
is well captured, with all the features observed in the laboratory
experiment reproduced in the virtual experimentation. 
Finally, we bring numerical evidence in support of the conjecture of
equivalence between our dynamic problem and the Monge-Kantorovich
equations by numerically solving classical Optimal Transport test
cases~\citep{prigozhin}. 
The numerical results are in well agreement with the results shown in
the literature and are obtained with a discretization scheme that is
highly attractive for its remarkable computational attributes of
simplicity and robustness.

\section*{Acknowledgments}

The authors are profoundly indebted with Giuseppe De Marco for his
continuous stimulation and his important contributions during the
development of this work. 

\bibliographystyle{abbrvnat}
\bibliography{Strings,biblio}

\begin{thebibliography}{20}
\providecommand{\natexlab}[1]{#1}
\providecommand{\url}[1]{\texttt{#1}}
\expandafter\ifx\csname urlstyle\endcsname\relax
  \providecommand{\doi}[1]{doi: #1}\else
  \providecommand{\doi}{doi: \begingroup \urlstyle{rm}\Url}\fi

\bibitem[Adamatzky(2010)]{physarum:machine}
A.~Adamatzky.
\newblock \emph{Physarum Machines}.
\newblock Computers from Slime Mould. World Scientific, 2010.

\bibitem[Ambrosio(2004)]{Ambrosio:2004}
L.~Ambrosio.
\newblock {Lecture Notes on Optimal Transport Problems}.
\newblock In \emph{Lecture Notes in Mathematics}, pages 1--52. Springer Berlin
  Heidelberg, Berlin, Heidelberg, 2004.

\bibitem[Ambrosio et~al.(2010)Ambrosio, Carlotto, and Massaccesi]{ambrosio_ln}
L.~Ambrosio, A.~Carlotto, and A.~Massaccesi.
\newblock \emph{Lecture Notes on Partial Differential Equations}.
\newblock {\tt http://cvgmt.sns.it/paper/1280/}, 2010.
\newblock URL \url{http://cvgmt.sns.it/paper/1280/}.

\bibitem[Barrett and Prigozhin(2007)]{prigozhin}
J.~W. Barrett and L.~Prigozhin.
\newblock A mixed formulation of the {M}onge-{K}antorovich equations.
\newblock \emph{Math. Model. Num. Anal.}, 41\penalty0 (6):\penalty0 1041--1060,
  2007.

\bibitem[Bochev and Lehoucq(2005)]{Bochev:2005}
P.~Bochev and R.~B. Lehoucq.
\newblock {On the Finite Element Solution of the Pure Neumann Problem}.
\newblock \emph{SIAM Rev.}, 47\penalty0 (1):\penalty0 50--66, Jan. 2005.

\bibitem[Bonifaci et~al.(2012)Bonifaci, Mehlhorn, and Varma]{bonifaci:physarum}
V.~Bonifaci, K.~Mehlhorn, and G.~Varma.
\newblock Physarum can compute shortest paths.
\newblock \emph{J Theor Biol}, 309:\penalty0 121--133, Sept. 2012.

\bibitem[Bouchitt{\'e} et~al.(1997)Bouchitt{\'e}, Buttazzo, and
  Seppecher]{shape}
G.~Bouchitt{\'e}, G.~Buttazzo, and P.~Seppecher.
\newblock Shape optimization solutions via {M}onge-{K}antorovich equation.
\newblock \emph{C. R. Acad. Sci. Paris S\'er. I Math}, 324\penalty0
  (10):\penalty0 1185--1191, 1997.

\bibitem[Brezzi and Fortin(1991)]{brezzi}
F.~Brezzi and M.~Fortin.
\newblock \emph{Mixed and Hybrid Finite Element Methods}.
\newblock Springer-Verlag, Berlin, 1991.

\bibitem[Buttazzo and Stepanov(2003)]{buttazzo}
G.~Buttazzo and E.~Stepanov.
\newblock On regularity of transport density in the monge--kantorovich problem.
\newblock \emph{SIAM J. Control Optim}, 42\penalty0 (3):\penalty0 1044--1055,
  2003.

\bibitem[Campanato(1965)]{Campanato:1965}
S.~Campanato.
\newblock Equazioni ellittiche del ii$^o$ ordine e spazi
  $\mathfrak{L}^{(2,\lambda)}$.
\newblock \emph{Annali di Matematica}, 69\penalty0 (1):\penalty0 321--381,
  1965.

\bibitem[Colombo and Mingione(2015)]{ColomboMingione:2015}
M.~Colombo and G.~Mingione.
\newblock Bounded minimisers of double phase variational integrals.
\newblock \emph{Arch. Rational Mech. Anal.}, 218\penalty0 (1):\penalty0
  219--273, Mar. 2015.

\bibitem[De~Pascale and Pratelli(2004)]{DePascale:2004}
L.~De~Pascale and A.~Pratelli.
\newblock Sharp summability for monge transport density via interpolation.
\newblock \emph{ESAIM: COCV}, 10\penalty0 (4):\penalty0 549--552, Oct. 2004.

\bibitem[Evans and Gangbo(1999)]{evans2}
L.~C. Evans and W.~Gangbo.
\newblock Differential equations methods for the {M}onge-{K}antorovich mass
  transfer problem.
\newblock \emph{Mem. Am. Math. Soc.}, 137\penalty0 (653):\penalty0 0--0, 1999.

\bibitem[Feldman and McCann(2002)]{Feldman:2002}
M.~Feldman and R.~J. McCann.
\newblock Uniqueness and transport density in monge's mass transportation
  problem.
\newblock \emph{Calc Var}, 15\penalty0 (1):\penalty0 81--113, 2002.

\bibitem[Giaquinta and Martinazzi(2013)]{giaquinta}
M.~Giaquinta and L.~Martinazzi.
\newblock \emph{{An Introduction to the Regularity Theory for Elliptic Systems,
  Harmonic Maps and Minimal Graphs}}.
\newblock Springer Science {\&} Business Media, Pisa, July 2013.

\bibitem[Morrey~Jr(1954)]{Morrey:1954}
C.~B. Morrey~Jr.
\newblock Second-order elliptic systems of differential equations.
\newblock In \emph{Contributions to the theory of partial differential
  equations}, pages 101--159. Princeton University Press, Princeton, N. J.,
  1954.

\bibitem[Nakagaki et~al.(2000)Nakagaki, Yamada, and Toth]{nakagaki:maze}
T.~Nakagaki, H.~Yamada, and A.~Toth.
\newblock Maze-solving by an amoeboid organism.
\newblock \emph{Nature}, 407\penalty0 (6803):\penalty0 470--470, 2000.

\bibitem[Tero et~al.(2007)Tero, Kobayashi, and Nakagaki]{tero:model}
A.~Tero, R.~Kobayashi, and T.~Nakagaki.
\newblock A mathematical model for adaptive transport network in path finding
  by true slime mold.
\newblock \emph{J. Theor. Biol.}, 244\penalty0 (4):\penalty0 553--564, Feb.
  2007.

\bibitem[Tero et~al.(2010)Tero, Takagi, Saigusa, Ito, Bebber, Fricker, Yumiki,
  Kobayashi, and Nakagaki]{tokyo}
A.~Tero, S.~Takagi, T.~Saigusa, K.~Ito, D.~P. Bebber, M.~D. Fricker, K.~Yumiki,
  R.~Kobayashi, and T.~Nakagaki.
\newblock Rules for biologically inspired adaptive network design.
\newblock \emph{Science}, 327\penalty0 (5964):\penalty0 439--442, Jan. 2010.

\bibitem[Troianiello(1987)]{troi}
G.~Troianiello.
\newblock \emph{Elliptic Differential Equations and Obstacle Problems}.
\newblock University Series in Mathematics. Springer, 1987.

\end{thebibliography}

\end{document}